\documentclass[preprint,12pt]{elsarticle}

\usepackage{hyperref}
\usepackage{amsmath,amssymb,amsthm,mathrsfs,amsfonts,dsfont}
\usepackage{graphicx}
\usepackage{subfigure}

\allowdisplaybreaks[1]

\usepackage{url}

\journal{Journal of Computers and Mathematics with Applications \hspace{0.2cm}}

\begin{document}

\begin{frontmatter}



\title{PDE formulation of some SABR/LIBOR market models and its numerical solution with a sparse grid combination technique\tnoteref{mytitlenote}}
\tnotetext[mytitlenote]{Partially financed by Spanish Grant MTM2013-47800-C2-1-P and by Xunta de Galicia (Grant CN2014/044). First author has also been founded by a FPU Spanish Grant.}

\author[auth]{J. G. L\'opez-Salas}\ead{jose.lsalas@udc.es}

\author[auth]{C. V\'azquez\corref{cor1}}\ead{carlosv@udc.es}

\address[auth]{Department of Mathematics, Faculty of Informatics, Campus Elvi\~na s/n, 15071-A Coru\~na (Spain)}

\cortext[cor1]{Corresponding author. Tel.: +34 981167000; fax: +34 981167160.}

\begin{abstract}

SABR models have been used to incorporate stochastic volatility to LIBOR market models (LMM) in order to describe interest rate dynamics and price interest rate derivatives. From the numerical point of view, the pricing of derivatives with SABR/LIBOR market models (SABR/LMMs) is mainly carried out with Monte Carlo simulation. However, this approach could involve excessively long computational times. For first time in the literature, in the present paper we propose an alternative pricing based on partial differential equations (PDEs). Thus, we pose original PDE formulations associated to the SABR/LMMs proposed by Hagan \cite{haganSABRLIBOR}, Mercurio \& Morini \cite{mercurioMorini} and Rebonato \cite{rebonatoWhite}. Moreover, as the PDEs associated to these SABR/LMMs are high dimensional in space, traditional full grid methods (like standard finite differences or finite elements) are not able to price derivatives over more than three or four underlying interest rates. In order to overcome this curse of dimensionality, a sparse grid combination technique is proposed. A comparison between Monte Carlo simulation results and the ones obtained with the sparse grid technique illustrates the performance of the method.
\end{abstract}

\begin{keyword}


Stochastic volatility models\sep SABR/LIBOR market models \sep High dimensional PDEs \sep Sparse grids \sep Combination technique
\end{keyword}

\end{frontmatter}


\section{Introduction}

The LMM \cite{braceLMM,jamshidianLMM,miltersenLMM} has become the most popular interest rate model. The main reason is the agreement between this model and Black's formulas \cite{brigoMercurio}. The standard LIBOR market model considers constant volatilities for the forward rates. However, this is a very limited hypothesis since it is impossible to reproduce market volatility smiles.

Among the different stochastic volatility models offered in the literature, the SABR model proposed by Hagan, Kumar, Lesniewski and Woodward \cite{Hagan} in the year 2002 stands out for becoming the market standard to reproduce the price of European options. SABR is the acronym for Stochastic, Alpha, Beta and Rho, three of the four model parameters. The SABR model can not be used to price derivatives whose payoff depends on several forward rates. In fact, SABR model works in the terminal measure, under which both the forward rate and its volatility are martingales. This can always be done if we work with one forward rate in isolation at a time. Under this same measure, however, the process for another forward rate and for its volatility would not be driftless.

In order to allow LMM to fit market volatility smiles, different extensions of the LMM that incorporate the volatility smile by means of the SABR model were proposed. These models are known as SABR/LIBOR market models (SABR/LMMs). In this article we will deal with the models proposed by Hagan \cite{haganSABRLIBOR}, Mercurio and Morini \cite{mercurioMorini} and Rebonato \cite{rebonatoWhite}.

While Monte Carlo \cite{book:Glasserman2003} simulation remains the common choice for pricing interest rate derivatives within SABR/LMM setting, several difficulties motivate to address alternative approaches based on PDE formulations. The first issue is that the convergence of Monte Carlo methods, although it depends only very weakly on the dimension of the problem, is very slow. Indeed, if the standard deviation of the result using a single simulation is $\epsilon$ then the standard deviation of the error after $N$ simulations is $\epsilon/\sqrt{N}$. Therefore, to improve the accuracy of the solution by a factor of $10$, $100$ times as many simulations must be performed. The second drawback of Monte Carlo methods is the valuation of options with early-exercise, like in the case of the American options, due to the so-called ``Monte Carlo on Monte Carlo'' effect. Available Monte Carlo methods for American options are also quite costly, see \cite{longstaffSchwartz} for example. In contrast, the modification of the PDE to a linear complementarity problem is usually straightforward. Finally, the weakest point of Monte Carlo methods appears to be the computation of the sensitivities of the solution with respect to the underlyings, the so-called ``Greeks'', which are very used by traders, and are directly given by the partial derivatives of the PDE solution. Besides, path-dependent options, like barrier options, can be easily priced in the PDE context where only the boundary conditions need to be changed, in contrast to Monte Carlo methods, where Brownian bridge techniques \cite{gobetBarrier} must be applied.

In view of previous arguments, in the present paper we pose equivalent PDE formulations for the three above mentioned SABR/LMMs. As far as we now, this is the first time in the literature that these PDE formulations are posed. From the numerical point of view, one main difficulty in these PDE formulations lies in their high dimensionality in space-like variables. In order to cope with this so-called {\em curse of dimensionality} several methods are available in the literature, see \cite{gerstner,beylkin} for example, which can be put into three categories. The first group uses the Karhunen-Loeve transformation to reduce the stochastic differential equation to a lower dimensional equation, therefore this results in a lower dimensional PDE associated to the previously reduced SDE. The second category gathers those methods which try to reduce the dimension of the PDE itself, like for example dimension-wise decomposition algorithms. Finally, the third category groups the methods which reduce the complexity of the problem in the discretization layer, like for example the method of sparse grids, which we use in the present article.

The sparse grid method was originally developed by Smolyak \cite{smolyak}, who used it for numerical integration. It is mainly based on a hierarchical basis \cite{yserentantI,yserentantII}, a representation of a discrete function space which is equivalent to the conventional nodal basis, and a sparse tensor product construction. Zenger \cite{zenger} and Bungartz and Griebel \cite{bungartzGriebel} extended this idea and applied sparse grids to solve PDEs with finite elements, finite volumes and finite differences methods. Besides working directly in the hierarchical basis, the sparse grid can also be computed using the combination technique \cite{griebelZenger} by linearly combining solutions on traditional Cartesian grids with different mesh widths. This is the approach we follow in this article. Recently, this technique has been used for a financial application related to the pricing of basket options in \cite{hendricksheuerehrhardtgunter,kraus,resingerWittum}.

The paper is organized as follows. In Section \ref{sec:intRateDer} some basic concepts related to interest rate derivatives and the corresponding terminology and notation are introduced. In Section \ref{sec:pdeFormulations} we pose the PDE formulations for the SABR/LMMs. In Section \ref{finiteDifferenceMethod} we describe the use of a full grid finite differences scheme for the Mercurio and Morini model, the application of which is analogous for the other two SABR/LMMs. Numerical results show the limitations of the full grid method when the number of forward rates increases. Therefore, in Section \ref{sec:sparseGrids} we describe the sparse grid combination technique applied to the SABR/LMM and show numerical results that illustrate the behaviour of the method when the number of forward rates increases. For this purpose, a comparison with Monte Carlo simulation results is used when analytic expressions of the solution are not available, as it happens in most of the cases. Note that Monte Carlo techniques are the usual alternative to price with SABR/LMM.

\section{Interest rate derivatives. Caplets and swaptions} \label{sec:intRateDer}
This section provides a brief introduction to the interest rate derivatives we deal with in the present article, for a deeper study we refer the reader to \cite{brigoMercurio}. Interest rate derivatives consist of financial contracts that depend on some interest rates.

A zero coupon \textit{bond} with maturity at time $T$ is a contract that pays its holder one unit of currency at time $T$. The value of this product at time $t<T$ is denoted by $P(t,T)$, and is called the discount factor from time $T$ to time $t$. Note that $P(T,T) = 1$ for all $T$.

A \textit{tenor structure} is a set of ordered payment dates $\{T_i, \; i=0, \dots, N\}$, such that $$T_0<T_1<\ldots<T_{N-1}<T_N.$$ The time between the payment dates is denoted by $\tau_i = T_{i+1}-T_i$. In terms of the corresponding discount factor, a payment of $x$ units at time $T_i$ is worth $xP(t,T_i)$ at time $t<T_i$.

A \textit{forward} interest rate $F_i(t)$ is an interest rate we can contract in order to borrow or lend money during the future time period $[T_i,T_{i+1}]$, and can be expressed in terms of discount factors in the form: $$F_i(t) = F(t;T_i,T_{i+1}) = \dfrac{1}{\tau_i}\left(\dfrac{P(t,T_i)}{P(t,T_{i+1})}-1\right) \mbox{ where } t\leq T_i.$$

Conversely, the price of a bond at time $T_i$ that matures at $T_j$, $P(T_i,T_j)$, can be expressed in terms of forward LIBOR rates as follows:
$$P(T_i,T_j) = \displaystyle\prod_{k=i}^{j-1}\dfrac{1}{1+\tau_k F_k(T_i)}.$$

A \textit{caplet} is a European call option on a forward rate. If the caplet expires at time $T_{i+1}$, at that time we will receive the payoff $\tau_i (F_i(T_i)-K)^+$, so that its discounted payoff at time $t<T_{i+1}$ is given by
$$P(t,T_{i+1})\tau_i (F_i(T_i)-K)^+,$$ where $(\cdot)^+$ denotes the function $\max(\cdot,0)$ and $K$ is the strike of the contract, which is given by a fixed interest rate in the contract. If constant volatilities are assumed, the caplet pricing can be computed with Black's formula \cite{brigoMercurio}
$$P(t,T_{i+1}) \tau_i {\mbox{Bl}}\Big(K, F_i(t),\sigma_{Black}\sqrt{T_i - t}\Big),$$ where
\begin{equation}
{\mbox{Bl}}(K, F, \nu) = F \varPhi\big(d_1(K,F,\nu)\big) - K \varPhi\big(d_2(K,F,\nu)\big),
\label{eq:black}
\end{equation}
$$d_1(K,F,\nu) = \dfrac{\ln(F/K)+\nu^2/2}{\nu},$$
$$d_2(K,F,\nu) = \dfrac{\ln(F/K)-\nu^2/2}{\nu},$$
and $\sigma_{Black}$ is the constant volatility of the forward rate which can be retrieved from market quotes.

An interest rate \textit{swap} (IRS) is a contract to exchange interest payments at future fixed dates. At every time instant $T_{i+1}$ in a prescribed set of dates $T_{a+1},\ldots,T_{b}$ the contract holder pays a fixed interest rate $K$ and receives a floating interest rate at the LIBOR rate $F_i(T_i)$ fixed at time $T_i$. The discounted payoff at time $t<T_a$ of this swap can be expressed as
\begin{equation}
\mbox{IRS}(t;T_a,\ldots,T_b) = \displaystyle\sum_{i=a}^{b-1} P(t,T_{i+1})\tau_i (F_i(T_i)-K). \label{eq:swap}
\end{equation}

A European $T_a \times (T_b-T_a)$ \textit{swaption} is an option giving the right (and not the obligation) to enter a swap at the future time $T_a$, called the swaption maturity. The underlying swap length $T_b-T_a$ is referred as the tenor of the swaption. Therefore, the discounted swaption payoff to the current time $t$ is equal to $$P(t,T_a)\big(\mbox{IRS}(T_a;T_a,\ldots,T_b)\big)^+.$$

As indicated in the introduction, the main objective of this article is the pricing of the previously described interest rate derivatives in the framework of SABR/LMM which incorporates the stochastic volatility by means of PDE formulations. Note that in most cases there are no analytical formulas for the solution. Thus, the new models are posed and solved with suitable numerical methods that overcome the high dimension in space of the equations. Moreover, the proposed sparse grid technique is parallelized to make the approach computationally efficient. Although along the paper we concentrate on Mercurio and Morini SABR/LMM, the proposed methodology can be analogously applied to Rebonato and Hagan models.

\section{Derivation of the PDE from the stochastic processes} \label{sec:pdeFormulations}

In \cite{carlosVazquezSA_SABRLIBOR} the authors analyzed the three SABR/LIBOR market models proposed by Hagan, Mercurio \& Morini and Rebonato using Monte Carlo simulation and their implementation on GPUs in order to price several interest rate derivatives. They have concluded that the Mercurio \& Morini model is the one with the best performance: it is the easiest to calibrate, it achieves the best fit to swaption market prices and it results the fastest one in the pricing with Monte Carlo simulation. Taking into account these reasons, in the present article we mainly choose this model to pose the PDE formulation and develop its numerical solution with the proposed methods. Nevertheless, at the end of this section we also pose the PDEs for the models of Hagan and Rebonato.

In order to describe the SABR/LMM setting, we first consider a set of $N-1$ LIBOR forward rates $F_i$, $1\leq i \leq N-1$, $\mathbf{F}=(F_1,\ldots,F_{N-1})$ on the tenor structure $[T_0,T_1, \ldots, T_{N-1},T_N]$, the accruals being $\tau_i = T_{i+1}-T_i$. The Mercurio \& Morini model is defined by the following system of stochastic differential equations \cite{mercurioMorini}:
\begin{align}
 dF_i(t) &= \mu_{i}(t) F_i(t)^\beta dt + \alpha_i V(t) F_i(t)^\beta dW^{\mathcal{Q}}_i(t), \quad F_i(0) \mbox{ given}, \nonumber \\
 dV(t) &= \sigma V(t) dZ^{\mathcal{Q}}(t), \quad V(0) = \alpha, \label{eq:merMorSDE}
\end{align}
which are posed on a probability space $\{\Omega, \mathcal{F}, \mathcal{Q}\}$ with filtration $\{\mathcal{F}_t\}$, $t\in[T_0,T_N]$. In (\ref{eq:merMorSDE}) $\mu_i$ is the drift of the $i$-th forward rate, $\beta \in [0,1]$ is the local volatility coefficient, $\alpha_i$ is a deterministic (constant) instantaneous volatility coefficient, $W_i^{\mathcal{Q}}$ are standard Brownian motions under the risk neutral measure $\mathcal{Q}$, $\rho$ is the correlation matrix between the forward rates, i.e.
$$<dW_i^\mathcal{Q}(t),dW_j^\mathcal{Q}(t)> = \rho_{ij}dt, \quad \forall i,j \in\{1,\ldots,N-1\},$$
$V$ is the stochastic volatility of the forward rates, $dZ^\mathcal{Q}$ is a standard Brownian motion correlated with the Brownian motions of the forward rates and $\phi$ is the correlation vector between the forward rates and the stochastic volatility, i.e.
$$<dW_i^\mathcal{Q}(t),dZ^\mathcal{Q}(t)>=\phi_i dt, \quad \forall i \in\{1, \ldots,N-1\}.$$
Due to the fact that the volatility process is lognormal, one can set the initial value of the volatility equal to one, i.e. $\alpha = 1$ with no loss of generality, since any different initial value can be embedded in the model by adjusting the deterministic coefficients $\alpha_i$. This is the choice we adopt in the following.

The drifts of the forward rates are determined by the chosen numeraire. Under the terminal probability measure $\mathcal{Q}^{T_N}$ associated with choosing the bond $P(t,T_N)$ as numeraire, the drifts of the forwards rates are given by
  $$\mu_{i}(t) = \left\{ \begin{array}{lcc}
             -\alpha_i V(t)^2 \displaystyle\sum_{j=i+1}^{N-1}\dfrac{\tau_j F_j(t)^\beta}{1+\tau_j F_j(t)}\rho_{ij}\alpha_j &   \mbox{if}  & i < N-1, \\
             \\ 0 &  \mbox{if}  & i = N-1.
             \end{array}
   \right.$$
In order to price $T_a\times(T_b-T_a)$ swaptions we will consider the probability measure $\mathcal{Q}^{T_a}$ associated with choosing the bond $P(t,T_a)$ as numeraire. In this case the drifts of the forward rates are given by
  $$\mu_{i}(t) = \left\{ \begin{array}{lcc}
               0 &  \mbox{if}  & i = a, \\
             \alpha_i V(t)^2 \displaystyle\sum_{j=a+1}^{i}\dfrac{\tau_j F_j(t)^\beta}{1+\tau_j F_j(t)}\rho_{ij}\alpha_j &   \mbox{if}  & i > a.
             \end{array}
   \right.$$

Our model for the correlation structure is taken from Rebonato \cite{rebonatoLMM}, who suggests the time independent function
\begin{equation}\rho_{ij} = \mbox{e}^{-\lambda|T_i-T_j|}.\label{eq:merMorCorr}\end{equation}
This function reflects the fact that the correlation increases as the time between the forward rates expiry decreases, so that two consecutive forward rates influence each other more than a forward rate in many years time.

A European option is characterized by its payoff function $G$, which determines the amount $G(T,\mathbf{F}(T),V(T))$ its holder receives at time $t=T$. The arbitrage-free value of the option relative to a numeraire $\mathcal{N}$ is then given by
\begin{equation}u(t,\mathbf{F}(t),V(t))=\dfrac{U(t,\mathbf{F}(t),V(t))}{N(t)} = \mathds{E}^\mathcal{Q}\left( \dfrac{G(T,\mathbf{F}(T),V(T))}{\mathcal{N}(T)} \Biggm| \mathcal{F}_t \right). \label{eq:merMorNum}\end{equation}
Closed-form solutions based on \eqref{eq:merMorNum} are rarely available due to the multi-asset feature of most LIBOR derivatives. In the next paragraphs we sketch the derivation of the PDE formulation associated to the Mercurio \& Morini model.

By using It\^o's formula, see \cite{steven} for example, the stochastic differential equation for $u$ is given by
\begin{align} du(t,\mathbf{F}(t),V(t)) \mbox{ } =& \mbox{ } \dfrac{\partial u}{\partial t}dt + \displaystyle\sum_{i=1}^{N-1} \dfrac{\partial u}{\partial F_i(t)} dF_i(t) + \dfrac{\partial u}{\partial V(t)} dV(t) +  \nonumber \\
 & \mbox{ } \dfrac{1}{2}\displaystyle\sum_{i,j=1}^{N-1} \dfrac{\partial^2 u}{\partial F_i(t) \partial F_j(t)} dF_i(t) dF_j(t) + \dfrac{1}{2} \dfrac{\partial^2 u}{\partial V(t)^2} (dV(t))^2 + \nonumber \\
 & \mbox{ } \displaystyle\sum_{i=1}^{N-1}\dfrac{\partial^2 u}{\partial F_i(t) \partial V(t)} dF_i(t) dV(t), \label{eq:merMorIto}
\end{align}
with box algebra \cite{mcLeish}: \\
\begin{center}
\begin{tabular}{c|cccc}
 & $dt$ & $dW^\mathcal{Q}_i$ & $dW^\mathcal{Q}_j$ & $dZ^\mathcal{Q}$ \\
\hline
 $dt$   & $0$ & $0$           & $0$           & $0$ \\
 $dW^\mathcal{Q}_i$ & $0$ & $dt$          & $\rho_{ij}dt$ & $\phi_i dt$ \\
 $dW^\mathcal{Q}_j$ & $0$ & $\rho_{ij}dt$ & $dt$          & $\phi_j dt$ \\
 $dZ^\mathcal{Q}$   & $0$ & $\phi_i dt$   & $\phi_j dt$   & $dt$ \\
\end{tabular}.
\end{center}
The interpretation of the box algebra is the following. In an expansion to terms of order $dt$, as $dt \rightarrow 0$ higher order terms such as $(dt)^j$ are all negligible for $j>0$. For example, $(dt)^2$ is of order $0$ as $dt\rightarrow0$, which is denoted as $(dt) (dt) \sim 0 $. Similarly, cross terms such as $(dt)(dW_i^\mathcal{Q})$ are negligible because the increment $dW^\mathcal{Q}_i$ is normally distributed with mean $0$ and standard deviation $(dt)^{1/2}$ and so $(dt)(dW^\mathcal{Q}_i)$ has standard deviation $(dt)^{3/2}$ which tends to $0$ as $dt\rightarrow0$.

Substituting equations \eqref{eq:merMorSDE} in \eqref{eq:merMorIto} and using the box algebra, we get
\begin{align}
du(t,\mathbf{F}(t),V(t)) = & \Biggm( \dfrac{\partial u}{\partial t} + \displaystyle \sum_{i=1}^{N-1} \mu_i(t) F_i(t)^\beta \dfrac{\partial u}{\partial F_i(t)} + \nonumber \\
 & \dfrac{1}{2} \displaystyle\sum_{i,j=1}^{N-1} \alpha_i \alpha_j V(t)^2 F_i(t)^\beta F_j(t)^\beta \rho_{ij} \dfrac{\partial^2 u}{\partial F_i(t) \partial F_j(t)} + \nonumber \\
 & \dfrac{1}{2}\sigma^2 V(t)^2 \dfrac{\partial^2 u}{\partial V(t)^2} + \displaystyle \sum_{i=1}^{N-1} \sigma V(t)^2 \alpha_i F_i(t)^\beta \phi_i \dfrac{\partial^2 u}{\partial F_i(t) \partial V(t)}\Biggm) dt + \nonumber \\
 & \displaystyle\sum_{i=1}^{N-1}\alpha_i V(t) F_i(t)^\beta \dfrac{\partial u}{\partial F_i(t)} dW_i^\mathcal{Q} + \sigma V(t) \dfrac{\partial u}{\partial V(t)} dZ^\mathcal{Q}. \label{eq:merMorProcess}
\end{align}

In order to comply with the no-arbitrage conditions and \eqref{eq:merMorNum}, the process $du(t,\mathbf{F},V)$ has to be martingale under the measure $\mathcal{Q}$. Thus, to satisfy this requirement, the drift term $dt$ in \eqref{eq:merMorProcess} must be equal to zero. The same result could be directly obtained by applying Feynman-Kac theorem, see \cite{steven,brigoMercurio}. The final parabolic PDE takes the following form:
\begin{eqnarray}
 \dfrac{\partial u}{\partial t} + \dfrac{1}{2} \sigma^2 V^2 \dfrac{\partial^2 u}{\partial V^2} + \dfrac{1}{2} V^2 \displaystyle \sum_{i,j=1}^{N-1} \rho_{ij} \alpha_i \alpha_j F_i^\beta F_j^\beta \dfrac{\partial^2 u}{\partial F_i \partial F_j} +  & & \nonumber \\  \sigma V^2 \displaystyle \sum_{i=1}^{N-1} \phi_{i} \alpha_i F_i^{\beta} \dfrac{\partial^2 u}{\partial F_i \partial V} + \displaystyle\sum_{i=1}^{N-1}\mu_i(t) F_i^\beta \dfrac{\partial u}{\partial F_i}  &=& 0, \label{eq:merMorPDE}
\end{eqnarray}
with the terminal condition given by the derivative payoff, $$u(T,\mathbf{F},V) = g(T,\mathbf{F},V),$$ on $\mathbb{R}^{N-1}\times\mathbb{R}$. For simplicity of notation, we have used the relative payoff $g(\cdot) = \dfrac{G(\cdot)}{\mathcal{N}(T)}$. The derivative price at time $t<T$ is given by $\mathcal{N}(t)u(t,\mathbf{F}(t),V(t))$.

Analytic solutions for \eqref{eq:merMorPDE} can be only found for suitable simple specifications of the functionals forms of the PDE and for straightforward boundary conditions (e.g. simple caplets without stochastic volatility, i.e. $\sigma = 0$, see Section \ref{numericaResutsDF}).

Finally, we are going to present the PDE for Hagan model, which is defined by the following system of stochastic differential equations \cite{haganSABRLIBOR}:
\begin{align}
 dF_i(t) &= \mu^{F_i}(t) F_i(t)^{\beta_i} dt + V_i(t) F_i(t)^{\beta_i} dW^{\mathcal{Q}}_i(t), \quad F_i(0) \mbox{ given}, \nonumber \\
 dV_i(t) &= \mu^{V_i}(t) V_i(t) dt + \sigma_i V_i(t) dZ_i^{\mathcal{Q}}(t), \quad V_i(0) = \alpha_i, \label{eq:haganSDE}
\end{align}
with the associated correlations denoted by
\begin{align*}
 <dW_i^\mathcal{Q}(t),dW_j^\mathcal{Q}(t)> &= \rho_{ij}dt,  \\
 <dW_i^\mathcal{Q}(t),dZ_j^\mathcal{Q}(t)> &= \phi_{ij}dt,  \\
 <dZ_i^\mathcal{Q}(t),dZ_j^\mathcal{Q}(t)> &= \theta_{ij}dt.
\end{align*}
The PDE for this model is obtained in the same way as previously with the Mercurio \& Morini model, thus obtaining:
\begin{align}
 \dfrac{\partial u}{\partial t} &+ \dfrac{1}{2} \displaystyle \sum_{i,j=1}^{N-1} \theta_{ij}\sigma_i V_i \sigma_j V_j \dfrac{\partial^2 u}{\partial V_i \partial V_j} + \dfrac{1}{2} \displaystyle \sum_{i,j=1}^{N-1} \rho_{ij} V_i F_i^{\beta_i} V_j F_j^{\beta_j} \dfrac{\partial^2 u}{\partial F_i \partial F_j} +  \nonumber \\
 &  \displaystyle \sum_{i,j=1}^{N-1} \phi_{ij} V_i F_i^{\beta_i} \sigma_j V_j \dfrac{\partial^2 u}{\partial F_i \partial V_j} + \displaystyle\sum_{i=1}^{N-1}\mu^{F_i}(t) F_i^{\beta_i} \dfrac{\partial u}{\partial F_i} + \displaystyle\sum_{i=1}^{N-1}\mu^{V_i}(t) V_i \dfrac{\partial u}{\partial V_i}= 0. \label{eq:haganPDE}
\end{align}

Rebonato model \cite{rebonatoWhite} is analogous to Hagan one, therefore its PDE will be also quite similar to \eqref{eq:haganPDE}.

\section{Finite Differences Method with full grids} \label{finiteDifferenceMethod}

Hereafter, as we have motivated in the previous section, we are going to just focus on the PDE \eqref{eq:merMorPDE} of the Mercurio \& Morini model. This backward parabolic PDE must be supplemented with a terminal condition, which describes the value of the variable $u$ at the final time $T$. Moreover, appropriate boundary conditions are required, which prescribe how the function $u$, or its derivatives, behave at the boundaries of the necessarily bounded computational domain.

We are going to define a $(N+1)$-dimensional mesh with the time sampled from today (time $0$) to the final expiry of the option (time $T$) at $M+1$ points uniformly spaced by the time step $\Delta t = \dfrac{T}{M}$.

The variables representing the forward rates $\textbf{F} = (F_1, \ldots, F_{N-1})$ and their stochastic volatility $V$, often referred as the ``space variables'' will be sampled at $M_i + 1$ ($i=1,\ldots,N-1$) and $S+1$ points spaced by $h_i = \dfrac{F_i^{max}-F_i^{min}}{M_i}$ and $h_v = \dfrac{V^{max}-V^{min}}{S}$, respectively.

Notice that while the choice of the range of the time variable is totally unambiguous, $[0,T]$, an a priori choice must be made about which values of the space variables are too high or too low to be of interest, so far we will denote them by $[F_i^{min},F_i^{max}]$ and $[V^{min},V^{max}]$. Selecting boundary values such that the option of interest is too deeply in or out-of-the money is a common and reasonable choice.

For a given mesh, each point is uniquely determined by the time level $m$ ($m=0,\ldots,M$), the index vector of the $N-1$ forward rates $\mathbf{f} = (f_1,\ldots,f_i, \ldots,$ $f_{N-1})$ ($f_i=0,\ldots,M_i$) and the stochastic volatility level $v$ ($v=0,\ldots,S$). We seek approximations of the solution at these mesh points, which will be denoted by $$U^m_{\mathbf{f},v} \approx u(m\Delta t, (f_i h_i)_{1\leq i \leq N-1}, v h_v).$$

It is natural for this PDE to be solved backwards in time. We approximate the time derivative by the time-forward approximation $$\dfrac{\partial u}{\partial t} \Biggm|_{t=m\Delta t, \mathbf{F}=(f_i h_i)_{1\leq i \leq N-1}, V = v h_v} = \dfrac{\partial u}{\partial t} \Biggm|_{m,\textbf{f},v} \approx \dfrac{U^{m+1}_{\mathbf{f},v} - U^{m}_{\mathbf{f},v}}{\Delta t}.$$

For the space derivatives we have chosen second-order approximations. We will write $\mathbf{f}_{i\pm1}$ to mean the forward rates index vector $(f_1,\ldots,f_i\pm1,\ldots,f_{N-1})$ which corresponds to the forward rates point $(f_1 h_1, \ldots, (f_i \pm 1) h_i, \ldots, f_{N-1} h_{N-1})$.

\noindent The first derivatives are approximated by central differences:
$$\dfrac{\partial u}{\partial F_i}\Biggm|_{m,\textbf{f},v} \approx \dfrac{U^{m}_{\mathbf{f}_{i+1},v} - U^{m}_{\mathbf{f}_{i-1},v}}{2h_i}.$$

\noindent The second derivatives are approximated by:
\begin{itemize}
 \item $\dfrac{\partial^2 u}{\partial F_i^2}\Biggm|_{m,\textbf{f},v} \approx \dfrac{U^{m}_{\mathbf{f}_{i+1},v} - 2 U^{m}_{\mathbf{f}_{i},v} + U^{m}_{\mathbf{f}_{i-1},v}}{h_i^2},$
 \item $\dfrac{\partial^2 u}{\partial V^2}\Biggm|_{m,\textbf{f},v} \approx \dfrac{U^{m}_{\mathbf{f},v+1} - 2 U^{m}_{\mathbf{f},v} + U^{m}_{\mathbf{f},v-1}}{h_v^2}.$
\end{itemize}

\noindent The cross derivatives terms are approximated by:
\begin{itemize}
 \item For $i\neq j$, $\dfrac{\partial^2 u}{\partial F_i \partial F_j}\Biggm|_{m,\textbf{f},v} \approx \dfrac{U^{m}_{\mathbf{f}_{i+1,j+1},v} + U^{m}_{\mathbf{f}_{i-1,j-1},v} - U^{m}_{\mathbf{f}_{i+1,j-1},v} - U^{m}_{\mathbf{f}_{i-1,j+1},v}}{4h_ih_j},$
 \item $\dfrac{\partial^2 u}{\partial F_i \partial V}\Biggm|_{m,\textbf{f},v} \approx \dfrac{U^{m}_{\mathbf{f}_{i+1},v+1} + U^{m}_{\mathbf{f}_{i-1},v-1} - U^{m}_{\mathbf{f}_{i+1},v-1} - U^{m}_{\mathbf{f}_{i-1},v+1}}{4h_ih_v}.$
\end{itemize}

\noindent The finite differences solution under the so-called $\theta$-scheme is:
$$ \dfrac{U^{m+1}_{\mathbf{f},v} - U^{m}_{\mathbf{f},v}}{\Delta t} + \theta W^m_{\mathbf{f},v} + (1-\theta) W^{m+1}_{\mathbf{f},v}=0,$$ where $\theta \in [0,1]$ and $W^m_{\mathbf{f},v}$ is the discretization given by
\begin{align}
 W^m_{\textbf{f},v} = & \nonumber \dfrac{1}{2}\sigma^2 V^2 \dfrac{U^{m}_{\mathbf{f},v+1} - 2 U^{m}_{\mathbf{f},v} + U^{m}_{\mathbf{f},v-1}}{h_v^2} + \nonumber \\
 & \dfrac{1}{2}V^2 \displaystyle \mathop{\sum_{i,j=1}^{N-1}}_{i\neq j} \rho_{ij} \alpha_i \alpha_j F_i^\beta F_j^\beta \dfrac{U^{m}_{\mathbf{f}_{i+1,j+1},v} + U^{m}_{\mathbf{f}_{i-1,j-1},v} - U^{m}_{\mathbf{f}_{i+1,j-1},v} - U^{m}_{\mathbf{f}_{i-1,j+1},v}}{4h_ih_j} + \nonumber \\
 & \dfrac{1}{2} V^2 \displaystyle\sum_{i=1}^{N-1} \alpha_i^2 F_i^{2\beta} \dfrac{U^{m}_{\mathbf{f}_{i+1},v} - 2 U^{m}_{\mathbf{f}_{i},v} + U^{m}_{\mathbf{f}_{i-1},v}}{h_i^2} + \nonumber \\
 & \sigma V^2 \displaystyle \sum_{i=1}^{N-1}\phi_i \alpha_i F_i^\beta \dfrac{U^{m}_{\mathbf{f}_{i+1},v+1} + U^{m}_{\mathbf{f}_{i-1},v-1} - U^{m}_{\mathbf{f}_{i+1},v-1} - U^{m}_{\mathbf{f}_{i-1},v+1}}{4h_i h_v} +\nonumber \\
 & \sum_{i=1}^{N-1} \mu_i(m \Delta t) F_i^\beta \dfrac{U^{m}_{\mathbf{f}_{i+1},v} - U^{m}_{\mathbf{f}_{i-1},v}}{2h_i}, \label{eq:merMorThetaScheme}
\end{align}
and with terminal condition $U^M_{\mathbf{f},v} = g(T,\mathbf{F},V)$, where we have denoted $\mathbf{F} = (F_i = f_i h_i)_{1 \leq i \leq N-1}$ and $V = v  h_v$.

Three different $\theta$ values represent three canonical discretization schemes, $\theta=0$ is the explicit scheme, $\theta=1$ the fully implicit scheme and $\theta = 0.5$ the Crank-Nicolson scheme. The fully implicit discretization is the best method with respect to stability, whereas the Crank-Nicolson timestepping provides the best convergence rate. Although the explicit method is the simplest to implement, it has the disadvantage of not being unconditionally stable.

We shall first discriminate explicit and implicit parts as follows:
\begin{equation}
 \dfrac{U_{\mathbf{f},v}^m}{\Delta t} - \theta W_{\mathbf{f},v}^m = \dfrac{U_{\mathbf{f},v}^{m+1}}{\Delta t} + (1-\theta) W_{\mathbf{f},v}^{m+1} .\label{eq:merMorSystem}
\end{equation}

As a result of such discretization we arrive to the linear system of equations $\mathbf{A}\mathbf{x} = \mathbf{b}$, where $\mathbf{A}$ is the band matrix of known coefficients, $\mathbf{x}$ is the vector of the unknown solutions $U^{m}_{\mathbf{f},v}$ and $\mathbf{b}$ is the vector of known values corresponding to the right-hand side of \eqref{eq:merMorSystem}.

Equation \eqref{eq:merMorSystem} can be rewritten as:
\begin{align}
 & d \theta U_{\mathbf{f},v-1}^m + d \theta U_{\mathbf{f},v+1}^m + \displaystyle \sum_{i=1}^{N-1} (b_i - r_i) \theta U_{\mathbf{f}_{i-1},v}^m + \displaystyle \sum_{i=1}^{N-1} (b_i + r_i) \theta U_{\mathbf{f}_{i+1},v}^m + \nonumber \\
 & \displaystyle \sum_{i=1}^{N-1} \big(a_i \theta U_{\mathbf{f}_{i-1},v-1}^m + a_i \theta U_{\mathbf{f}_{i+1},v+1}^m - a_i \theta U_{\mathbf{f}_{i-1},v+1}^m - a_i \theta U_{\mathbf{f}_{i+1},v-1}^m\big) + \nonumber \\
 & \displaystyle \sum_{ij \in C} \big( \psi_{ij} \theta U_{\mathbf{f}_{i-1,j-1},v}^m + \psi_{ij} \theta U_{\mathbf{f}_{i+1,j+1},v}^m - \psi_{ij} \theta U_{\mathbf{f}_{i-1,j+1},v}^m - \psi_{ij} \theta U_{\mathbf{f}_{i+1,j-1},v}^m \big) + \nonumber \\
 & \left( -1 - 2 d \theta - 2 \theta \sum_{i=1}^{N-1}b_i \right) U_{\mathbf{f},v}^m = \nonumber \\
 & - d \hat\theta U_{\mathbf{f},v-1}^{m+1} - d \hat\theta U_{\mathbf{f},v+1}^{m+1} - \displaystyle \sum_{i=1}^{N-1} (b_i - r_i) \hat\theta U_{\mathbf{f}_{i-1},v}^{m+1} - \displaystyle \sum_{i=1}^{N-1} (b_i + r_i) \hat\theta U_{\mathbf{f}_{i+1},v}^{m+1}  \nonumber \\
 & - \displaystyle \sum_{i=1}^{N-1} \big(a_i \hat\theta U_{\mathbf{f}_{i-1},v-1}^m + a_i \hat\theta U_{\mathbf{f}_{i+1},v+1}^{m+1} - a_i \hat\theta U_{\mathbf{f}_{i-1},v+1}^{m+1} - a_i \hat\theta U_{\mathbf{f}_{i+1},v-1}^{m+1}\big)  \nonumber \\
 & -\displaystyle \sum_{ij \in C} \big( \psi_{ij} \hat\theta U_{\mathbf{f}_{i-1,j-1},v}^{m+1} + \psi_{ij} \hat\theta U_{\mathbf{f}_{i+1,j+1},v}^{m+1} - \psi_{ij} \hat\theta U_{\mathbf{f}_{i-1,j+1},v}^{m+1} - \psi_{ij} \hat\theta U_{\mathbf{f}_{i+1,j-1},v}^{m+1} \big) + \nonumber \\
 & \left( -1 + 2 d \hat\theta + 2 \hat\theta \sum_{i=1}^{N-1}b_i \right) U_{\mathbf{f},v}^{m+1}, \label{eq:merMorDisc}
\end{align}
where $\hat\theta = (1-\theta)$, $C$ is the set containing the combinations of numbers $1,2,\ldots,N-1$ taken two at a time without repetition (the number of elements in $C$ is ${N-1 \choose 2} = 2^{-1}(N-1)(N-2)$) and the known coefficients $d$, $b_i$, $r_i$, $a_i$ and $\psi_{ij}$ are defined as
\begin{align*}
 & d = \dfrac{\Delta t\sigma^2V^2}{2h_v^2}, \quad b_i = \dfrac{\Delta t V^2 \alpha_i^2 F_i^{2\beta}}{2h_i^2}, \\
 & r_i = \dfrac{\Delta t \mu_i(t)F_i^\beta}{2h_i}, \quad a_i = \dfrac{\Delta t \sigma V^2 \phi_i \alpha_i F_i^\beta}{4h_ih_v}, \\
 & \psi_{ij} = \dfrac{\Delta t V^2 \rho_{ij} \alpha_i\alpha_j F_i^\beta F_j^\beta}{4h_ih_j}.
\end{align*}

\subsection{Boundary conditions} \label{boundaryConditions}

In order to specify boundary conditions, a combination of mathematical, financial and heuristic reasoning allows us to find consistent and acceptable ones. There are several possibilities, see \cite{duffy} for example.

 We assume that forward rates and their stochastic volatility are non negative and hence take values in the range zero to infinity. We first truncate the unbounded interval to a bounded one and then we must specify conditions at the new boundary. Thus we will consider the truncated domain $[F_i^{min}, F_i^{max}]\times[V^{min},V^{max}]$, with $F_i^{min} = 0$ and $V^{min}=0$.

For the forward rates we consider Dirichlet boundary conditions. Particularly, the terminal condition holds on the forward rates boundaries, i.e. $$U^m_{\{\mathbf{f} | \exists f_i = 0 \},v} = U^M_{\mathbf{f},v}, \quad \forall m=0,\ldots,M-1,$$ $$U^m_{\{\mathbf{f} | \exists f_i = M_i \},v} = U^M_{\mathbf{f},v}, \quad \forall m=0,\ldots,M-1.$$

At the stochastic volatility boundaries we consider the following conditions:
\begin{align}
 &\dfrac{\partial u}{\partial t} + \displaystyle \sum_{i=1}^{N-1} \mu_i(t) F_i^\beta \dfrac{\partial u}{\partial F_i} = 0, \quad V=0, \label{eq:boundCondition0} \\
 &\dfrac{\partial u}{\partial V} = 0, \quad V=V_{max}. \label{eq:boundConditionMax}
\end{align}
Thus, when $V=0$ we require that the PDE itself must be satisfied on this boundary. When $V$ approaches to infinity, the price of the derivative becomes independent of $V$. This is reflected by using Neumann conditions instead of the Dirichlet ones used for the forward rates boundaries.

At the boundary $V=0$, after discretizing the boundary condition \eqref{eq:boundCondition0} we obtain (note that the coefficients $d$, $b_i$, $a_i$ and $\psi_{ij}$ of equation \eqref{eq:merMorDisc} are zero):
\begin{align*}
 -\displaystyle\sum_{i=1}^{N-1}r_i \theta U^m_{\mathbf{f}_{i-1},0} + \displaystyle\sum_{i=1}^{N-1}r_i \theta U^m_{\mathbf{f}_{i+1},0} - U^m_{\mathbf{f},0} = \displaystyle\sum_{i=1}^{N-1}r_i \hat\theta U^{m+1}_{\mathbf{f}_{i-1},0} - \displaystyle\sum_{i=1}^{N-1}r_i \hat\theta U^{m+1}_{\mathbf{f}_{i+1},0} + U^{m+1}_{\mathbf{f},0}.
\end{align*}

For the boundary $V = V_{max}$ in order to maintain the second order accuracy in the discretization of the first derivative the ghost point method is considered. The ghost grid points $U_{\mathbf{f},S+1}$ are added. Then, the finite differences scheme of equation \eqref{eq:merMorDisc} can also be applied at the points $U_{\mathbf{f},S}$. However, we now have more unknowns than equations. The additional equations come from the central finite differences discretization of the Neumann boundary condition \eqref{eq:boundConditionMax}:
$$\dfrac{U_{\mathbf{f},S+1}-U_{\mathbf{f},S-1}}{2h_v} = 0,$$ which yields $U_{\mathbf{f},S+1} = U_{\mathbf{f},S-1}$. Inserting this into the finite differences equation at $V=V_{max}$ we achieve
\begin{align*}
 & \hat d \theta U^{m}_{\mathbf{f},S-1} + \displaystyle \sum_{i=1}^{N-1} (b_i-r_i)\theta U^{m}_{\mathbf{f}_{i-1},S} + \displaystyle \sum_{i=1}^{N-1} (b_i+r_i)\theta U^{m}_{\mathbf{f}_{i+1},S} + \nonumber \\
 & \displaystyle \sum_{ij \in C} \big( \psi_{ij} \theta U_{\mathbf{f}_{i-1,j-1},S}^m + \psi_{ij} \theta U_{\mathbf{f}_{i+1,j+1},S}^m - \psi_{ij} \theta U_{\mathbf{f}_{i-1,j+1},S}^m - \psi_{ij} \theta U_{\mathbf{f}_{i+1,j-1},S}^m \big) + \nonumber \\
 & \left( -1 - \hat d \theta - 2 \theta \sum_{i=1}^{N-1}b_i \right) U_{\mathbf{f},S}^{m} = \nonumber \\
 & -\hat d \hat\theta U^{m+1}_{\mathbf{f},S-1} - \displaystyle \sum_{i=1}^{N-1} (b_i-r_i)\hat\theta U^{m+1}_{\mathbf{f}_{i-1},S} - \displaystyle \sum_{i=1}^{N-1} (b_i+r_i)\hat\theta U^{m+1}_{\mathbf{f}_{i+1},S} + \nonumber \\
 & -\displaystyle \sum_{ij \in C} \big( \psi_{ij} \hat\theta U_{\mathbf{f}_{i-1,j-1},S}^{m+1} + \psi_{ij} \hat\theta U_{\mathbf{f}_{i+1,j+1},S}^{m+1} - \psi_{ij} \hat\theta U_{\mathbf{f}_{i-1,j+1},S}^{m+1} - \psi_{ij} \hat\theta U_{\mathbf{f}_{i+1,j-1},S}^{m+1} \big) + \nonumber \\
 & \left( -1 + \hat d \hat\theta + 2 \hat\theta \sum_{i=1}^{N-1}b_i \right) U_{\mathbf{f},S}^{m+1},
\end{align*}
where $\hat d = 2 d = \dfrac{\Delta t \sigma^2 V_{max}^2}{h_v^2}$.

\subsection{Numerical results} \label{numericaResutsDF}

It is not clear where to place $F_i^{max}$ and $V^{max}$. On one hand, it is advantageous to place them far away of the initial forward rates. This reduces the error of the artificial boundary conditions. On the other hand a large computational domain requires a large discretization width. This increases the error of the approximation of the derivatives. In our experiments we will consider $F_i^{max} = 0.1$ and $V^{max} = 3.5$.

Some specifications of the financial product are given in Table \ref{table:specification} and the employed market data, taken from \cite{blackham}, are shown in Table \ref{table:marketData}. We will consider $\lambda = 0.1$ in the model for the correlation structure \eqref{eq:merMorCorr}. Besides, the Crank-Nicolson scheme will be used in \eqref{eq:merMorSystem}. For solving the system \eqref{eq:merMorDisc} the Gauss-Seidel iterative solver has been employed using a tolerance of $
10^{-6}$.

The numerical experiments have been performed with the following hardware and software configurations: two recent multicore Intel Xeon CPUs E5-2620 v2 clocked at 2.10 GHz (6 cores per socket) with 62 GBytes of RAM, CentOS Linux, GNU C++ compiler 4.8.2.

\begin{table}[h]
\begin{center}
\begin{tabular}{|r|r|}
\hline
Currency & EUR \\
\hline
 Index & EURIBOR \\ %
\hline
 Day Count & e30/360 \\
\hline
 Strike & $5.5\%$ \\
\hline
\end{tabular}
\end{center}
\caption{Specification of the interest rate model.}
\label{table:specification}
\end{table}

\begin{table}[h]
\begin{center}
\begin{tabular}{|r|r|r|r|r|}
\hline
 & Start date & End date & LIBOR Rate (\%) & Volatility (\%) \\
\hline
\hline
 $T_0$ & $29$-$07$-$04$ & $29$-$07$-$05$ & $2.423306$ & $0$ \\
\hline
 $T_1$ & $29$-$07$-$05$ & $29$-$07$-$06$ & $3.281384$ & $24.73$ \\
\hline
 $T_2$ & $29$-$07$-$06$ & $29$-$07$-$07$ & $3.931690$ & $22.45$ \\
\hline
 $T_3$ & $29$-$07$-$07$ & $29$-$07$-$08$ & $4.364818$ & $19.36$ \\
\hline
 $T_4$ & $29$-$07$-$08$ & $29$-$07$-$09$ & $4.680236$ & $17.43$ \\
\hline
 $T_5$ & $29$-$07$-$09$ & $29$-$07$-$10$ & $4.933085$ & $16.15$ \\
\hline
 $T_6$ & $29$-$07$-$10$ & $29$-$07$-$11$ & $5.135066$ & $15.02$ \\
\hline
 $T_7$ & $29$-$07$-$11$ & $29$-$07$-$12$ & $5.273314$ & $14.24$ \\
\hline
 $T_8$ & $29$-$07$-$12$ & $29$-$07$-$13$ & $5.376115$ & $13.42$ \\
\hline
\end{tabular}
\end{center}
\caption{Market data used in pricing. Data taken from $27$th July 2004.}
\label{table:marketData}
\end{table}

We are going to value $T_a \times (T_b - T_a)$ European swaptions under the Mercurio \& Morini model.

First of all, the results when pricing a $1\times 1$ European swaption are discussed. The value of this particular swaption is the same as the price of the corresponding caplet, and so it only depends on $F_1$. Hence, in one dimension a closed form expression for the price of a European swaption can be found by using Black's formula \eqref{eq:black} and is given by
$$P(T_0,T_2) \tau_1 {\mbox{Bl}}(K, F_1(T_0),\sigma_{Black}\sqrt{T_1-T_0}).$$
This value is equal to $0.659096$ basis points (one basis point is one hundredth of one percent, $\frac{1\%}{100}=10^{-4}$). As Black's formula for caplets considers constant volatility $\sigma_{Black}$, in this first test the volatility of the volatility parameter of Mercurio \& Morini model is considered equal to zero, i.e., $\sigma = 0$, therefore a standard LIBOR market model is used. The solution was found on several levels and Table \ref{table:df2_0} shows the convergence of the model when using $256$ time steps. In all tables of this article, {\it Level} refers to the refinement level $n$, i.e., the mesh size is $h_i = 2^{-n} \cdot c_i$ in each coordinate direction, where $c_i$ denotes the computational domain length in direction $i$, which is $F_i^{max}$ in the case of the forward rates and $V^{max}$ in the case of the stochastic volatility. Besides, the solution and the error with respect to the exact solution are also shown in basis points. Additionally, the execution time is measured in seconds and the column labeled as {\it Grid points} shows the number of grid points employed in the full grid used by the finite differences method without taking into account the time coordinate. The time discretization error decreases as increasing the number of time steps. As an example, in Table \ref{table:df2_0_12} we show the achieved solutions and the associated errors using only $12$ time steps for the higher resolution levels ($n=9, \,10$), these results are to be compared with those of Table \ref{table:df2_0}. From the results in these two Tables \ref{table:df2_0} and \ref{table:df2_0_12}, we can see that not only the time step $\Delta t = 1/256$ is small enough but also the reported errors are dominated by the spatial discretization, at least up to level $10$. Therefore, in the following numerical tests we will consider $256$ time steps.

When the volatility of the volatility $\sigma$ of the model is non zero or when the length of the underlying swap of the swaption being considered is greater than one, no closed form solutions are available. However, an estimate can be obtained from Monte Carlo simulations. On Table \ref{table:d2_mc} Monte Carlo values for the $1 \times 1$ European swaption with $\sigma=0$ are shown for several numbers of paths ({\it \#Paths}). More details about Monte Carlo simulation of SABR/LMMs can be found in the article \cite{carlosVazquezSA_SABRLIBOR}.

\begin{table}[!hb]
\begin{center}
{\footnotesize
\begin{tabular}{|r|r|r|r|r|}
\hline
Level & Solution & Error & Time & Grid points\\
\hline
\hline
 $3$ & $2.529194$ & $1.870098$ & $0.01$ & $81$ \\ %
\hline
 $4$ & $1.205248$ & $0.546151$ & $0.03$ & $289$ \\ %
\hline
 $5$ & $0.802744$ & $0.143647$ & $0.21$ & $1089$ \\ %
\hline
 $6$ & $0.674394$ & $0.015297$ & $1.95$ & $4225$ \\ %
\hline
 $7$ & $0.663856$ & $0.004760$ & $24.42$ & $16641$ \\ %
\hline
 $8$ & $0.659937$ & $0.000840$ & $340.57$ & $66049$ \\ %
\hline
 $9$ & $0.659381$ & $0.000285$ & $4775.27$ & $263169$ \\ %
\hline
 $10$ & $0.659163$ & $0.000066$ & $65013.50$ & $1050625$ \\ %
\hline
\end{tabular}
}
\end{center}
\caption{Convergence of the full grid finite differences solution in basis points for $1$ LIBOR and stochastic volatility, $\sigma=0$, $V(0)=1$, $\beta=1$, $256$ time steps. Exact solution, $0.659096$ basis points.}
\label{table:df2_0}
\end{table}

\begin{table}[!hb]
\begin{center}
{\footnotesize
\begin{tabular}{|r|r|r|r|r|}
\hline
Level & Solution & Error & Time & Grid points\\
\hline
\hline
 $9$ & $0.658014$ & $0.001082$ & $1884.24$ & $263169$ \\ %
\hline
 $10$ & $0.657795$ & $0.001301$ & $23436.59$ & $1050625$ \\ %
\hline
\end{tabular}
}
\end{center}
\caption{Convergence of the full grid finite differences solution in basis points for $1$ LIBOR and stochastic volatility, $\sigma=0$, $V(0)=1$, $\beta=1$, $12$ time steps. Exact solution, $0.659096$ basis points.}
\label{table:df2_0_12}
\end{table}

\begin{table}[!hb]
\begin{center}
{\footnotesize
\begin{tabular}{|r|r|r|}
\hline
\#Paths & CI \\
\hline
\hline
 $10^5$ & $[0.608003,0.703098]$ \\ %
 $10^7$ & $[0.654232,0.663567]$ \\ %
 $10^9$ & $[0.658629,0.659562]$\\
\hline
\end{tabular}
}
\end{center}
\caption{$95$\% confidence intervals (CI) with Monte Carlo solution in basis points for 1 LIBOR and stochastic volatility, $\sigma=0$, $V(0)=1$, $\beta=1$, $256$ time steps. Exact solution, $0.659096$ basis points.}
\label{table:d2_mc}
\end{table}

In Table \ref{table:df2_1} the pricing of the $1 \times 1$ European swaption with $\sigma=0.3$  for different resolution levels $n$ are shown. In Table \ref{table:df3} the results for the $1 \times 2$ swaption are given. Note that with this numerical method it was not feasible to price the swaption past refinement level $n=8$ due to the huge number of required grid points. In Table \ref{table:df4} the results for the $1 \times 3$ swaption are given. Full grid pricing is only possible on low grid levels. It is not achievable to obtain a solution for a level greater than 6 in reasonable computational time due to the high number of grid points.

\begin{table}[!htb]
\begin{center}
{\footnotesize
\begin{tabular}{|r|r|r|r|}
\hline
Level & Solution & Time & Grid points\\
\hline
\hline
 $3$ & $3.440$ & $0.01$ & $81$ \\ %
\hline
 $4$ & $2.168$ & $0.03$ & $289$ \\ %
\hline
 $5$ & $1.800$ & $0.21$ & $1089$ \\ %
\hline
 $6$ & $1.678$ & $2.20$ & $4225$ \\ %
\hline
 $7$ & $1.670$ & $27.57$ & $16641$ \\ %
\hline
 $8$ & $1.667$ & $376.91$ & $66049$ \\ %
\hline
 $9$ & $1.665$ & $5206.93$ & $263169$ \\ %
\hline
 $10$ & $1.663$ & $70492.91$ & $1050625$ \\ %
\hline
\end{tabular}
}
\end{center}
\caption{Convergence of the full grid finite differences solution in basis points for $1$ LIBOR and stochastic volatility, $\sigma=0.3$, $\phi_i = 0.4$, $V(0)=1$, $\beta=1$, $256$ time steps. $95$\% confidence interval with Monte Carlo simulation using $10^7$ paths, $[1.652, 1.672]$ in basis points.}
\label{table:df2_1}
\end{table}

\begin{table}[!htb]
\begin{center}
{\footnotesize
\begin{tabular}{|r|r|r|r|}
\hline
Level & Solution & Time & Grid points\\
\hline
\hline
 $3$ & $7.162$ & $0.14$ & $729$ \\ %
\hline
 $4$ & $5.565$ & $1.84$ & $4913$ \\ %
\hline
 $5$ & $5.003$ & $34.41$ & $35937$ \\ %
\hline
 $6$ & $4.865$ & $806.02$ & $274625$ \\ %
\hline
 $7$ & $4.846$ & $21903.33$ & $2146689$ \\ %
\hline
 $8$ & $4.824$ & $611725.64$ & $16974593$ \\ %
\hline
\end{tabular}
}
\end{center}
\caption{Convergence of the full grid finite differences solution in basis points for $2$ LIBORs and stochastic volatility, $\sigma=0.3$, $\phi_i=0.4$, $V(0)=1$, $\beta=1$, $256$ time steps. $95$\% confidence interval with Monte Carlo simulation using $10^7$ paths, $[4.800,4.844]$ in basis points.}
\label{table:df3}
\end{table}

\begin{table}[!htb]
\begin{center}
{\footnotesize
\begin{tabular}{|r|r|r|r|}
\hline
Level & Solution & Time & Grid points\\
\hline
\hline
 $3$ & $11.702$ & $2.17$ & $6561$ \\ %
\hline
 $4$ & $9.497$ & $73.90$ & $83521$ \\ %
\hline
 $5$ & $8.892$ & $3033.58$ & $1185921$ \\ %
\hline
 $6$ & $8.771$ & $157152.75$ & $17850625$ \\ %
\hline
\end{tabular}
}
\end{center}
\caption{Convergence of the full grid finite differences solution in basis points for $3$ LIBORs and stochastic volatility, $\sigma=0.3$, $\phi_i=0.4$, $V(0)=1$, $\beta=1$, $256$ time steps. $95$\% confidence interval with Monte Carlo simulation using $10^7$ paths, $[8.635,8.700]$ in basis points.}
\label{table:df4}
\end{table}

Theoretically, it is possible to solve the discrete system \eqref{eq:merMorDisc} for a general number of dimensions. However, in computational science, a major problem occurs when the number of dimensions increases. A natural way to diminish the discretization error is to decrease the mesh step in each coordinate direction. However, then the number of grid points in the resulting full grid grows exponentially with the dimension, i.e. the size of the discrete solution increases drastically. This is called the \textit{curse of dimensionality} \cite{bellmann}. Therefore, this procedure of improving the accuracy by decreasing the mesh step is mainly bounded by two factors, the storage and the computational complexity. Due to these limitations, using a full grid discretization method which achieves sufficiently accurate approximations is only possible for problems with up to three or four dimensions, even on the most powerful machines presently available \cite{bungartzGriebel}.

\section{Sparse grids and the combination technique} \label{sec:sparseGrids}
Two approaches to try to overcome the curse of dimensionality are increasing the order of accuracy of the applied numerical approximation scheme or reducing the dimension of the problem by choosing suitable coordinates. Both approaches are not always possible for every option pricing problem. In this article we will take advantage of the sparse grid combination technique first introduced by Zenger and co-workers \cite{griebelZenger} in order to try to overcome the curse of dimensionality and allow to use the PDE formulation of SABR/LMM for the pricing problem we are dealing with. The combination technique replicates the structure of a so-called sparse grid by linearly combining solutions on coarser grids of the same dimensionality. This technique reduces the computational effort and the storage space involved with the mentioned traditional finite differences discretization methods. The number of sub-problems to solve increases, while the computational time per problem decreases drastically. This method can be implemented in parallel as each sub-grid is independent of the others. In the next two subsections we give a brief introduction to sparse grids and the combination technique. For a detailed discussion we refer to \cite{bungartzGriebel}.

\subsection{Sparse grids}
First, we introduce some notations and definitions. Let $\mathbf{l} = (l_1,l_2,\ldots,l_d) \in \mathds{N}_0^d$ denote a $d$-dimensional multi-index. Let $|\mathbf{l}|_1$ and $|\mathbf{l}|_\infty$ denote the discrete $L_1-$norm and $L_\infty-$norm of the multi-index $\mathbf{l}$, respectively, that are defined as $$|\mathbf{l}|_1 = \sum_{k=1}^d l_k \quad \mbox{and} \quad |\mathbf{l}|_\infty = \max_{1\leq k \leq d} l_k.$$ We define the anisotropic grid $\Omega_\mathbf{l}$ with mesh size $\mathbf{h} = (h_1,h_2,\ldots,h_d) = (2^{-l_1} c_1,$ $2^{-l_2} c_2, \ldots, 2^{-l_d} c_d)$ with multi-index $\mathbf{l}$ and grid length $\mathbf{c} = (c_1,c_2,\ldots,c_d)$.

\begin{figure}[!htb]
\begin{center}
\includegraphics[scale=1.1]{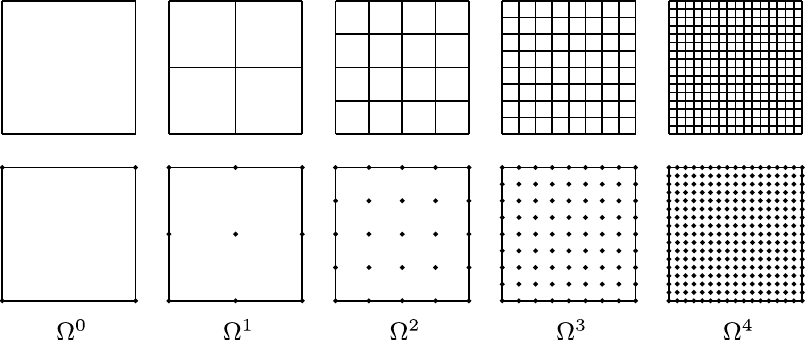}
\caption{Two-dimensional full grid hierarchy up to level $n=4$.}
\label{fig:merMorfullGrid}
\end{center}
\end{figure}

Then, the full grid at refinement level $n\in\mathds{N}$ and mesh size $h_i = 2^{-n} \cdot c_i$ for all $i$ can be defined via the sequence of subgrids $$\Omega^{n} = \Omega_{(n,\ldots,n)} = \bigcup\limits_{|\mathbf{l}|_\infty \leq n} \Omega_\mathbf{l}.$$
Figure \ref{fig:merMorfullGrid} visualizes two dimensional full grids for levels $n=0,\ldots,4$.

\noindent The number of grid points in each coordinate direction of the full grid is $2^n + 1$ and therefore the number of grid nodes in the full grid increases with $O(2^{n \cdot d})$, i.e. grows exponentially with the dimensionality $d$ of the problem.

\begin{figure}[!htb]
\begin{center}
\includegraphics[height=18cm]{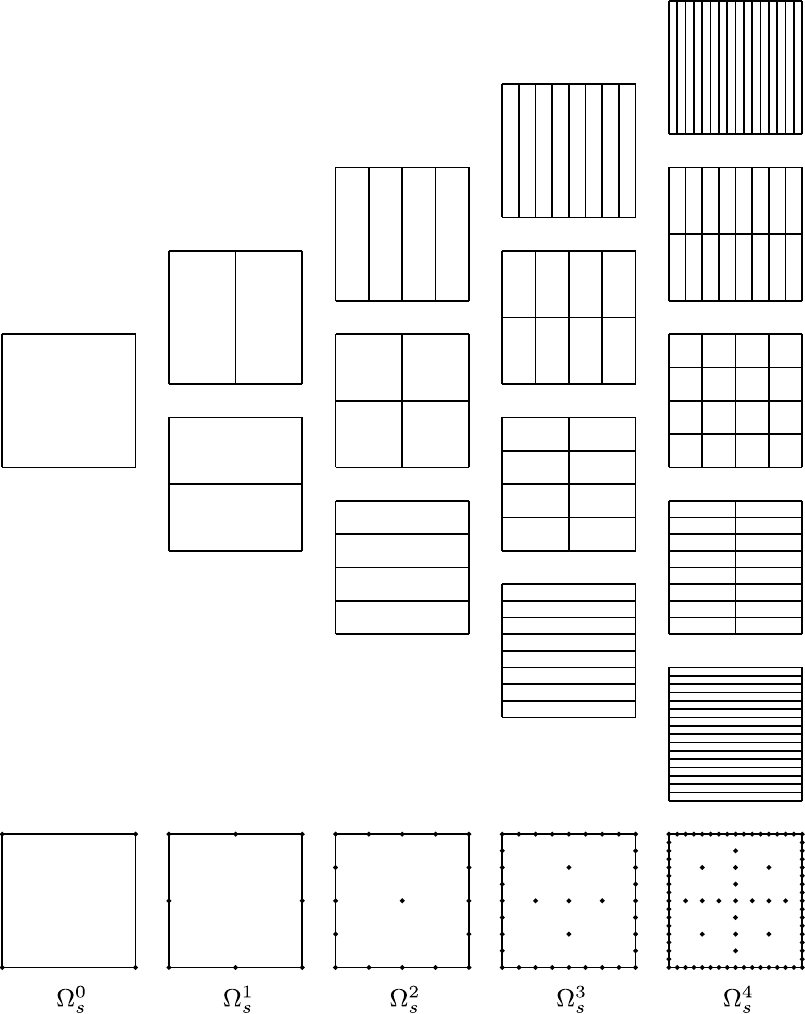}
\caption{Two-dimensional sparse grid hierarchy up to level $n=4$.}
\label{fig:merMorSparseGrid}
\end{center}
\end{figure}

The sparse grid $\Omega_s^n$ at refinement level $n$ consists of all anisotropic Cartesian grids $\Omega_{\mathbf{l}}$, where the total sum of all refinement factors $l_k$ in each coordinate direction equals the resolution $n$. Then, the sparse grid $\Omega_s^n$ is given by $$\Omega_s^n = \bigcup\limits_{|\mathbf{l}|_1\leq n} \Omega_\mathbf{l} = \bigcup\limits_{|\mathbf{l}|_1= n} \Omega_\mathbf{l}.$$
Figure \ref{fig:merMorSparseGrid} shows the two-dimensional grid hierarchy for levels $n=0,\ldots,4$.

The total number of nodes in the grid $\Omega_{\mathbf{l}}$ is $\displaystyle\prod_{k=1}^{d} (2^{l_k}+1) = O(2^{|\mathbf{l}|_1}) = O(2^n)$. In addition, there exist exactly $n+d-1 \choose d-1$ grids $\Omega_\mathbf{l}$ with $|\mathbf{l}|_1 = n$,
\begin{align*}{n+d-1 \choose d-1} &= \dfrac{(n+d-1)!}{(d-1)!n!} = \dfrac{(n+d-1)\cdot\ldots\cdot(n+1)n!}{(d-1)!n!}   \\
& = \dfrac{n+(d-1)}{d-1} \cdot \dfrac{n+(d-2)}{d-2} \cdot \ldots \cdot \dfrac{n+(d-(d-1))}{d-(d-1)} \\
& = \left(1 + \dfrac{n}{d-1} \right) \cdot \left(1 + \dfrac{n}{d-2} \right) \cdot \ldots \cdot \left(1 + \dfrac{n}{2}\right) \cdot \left(1 + \dfrac{n}{1} \right) \\
&\leq (1+n)^{d-1} = O(n^{d-1}).
\end{align*}
Thus, the total number of grid points of the sparse grid $\Omega_s^n$ grows according to
\begin{equation} {n+d-1 \choose d-1} \cdot \displaystyle\prod_{k=1}^{d} (2^{l_k}+1) = O(n^{d-1}) O(2^n) = O(n^{d-1} 2^n), \label{eq:merMorNumOfSGP}\end{equation}
which is far less the size of the corresponding full grid with $O(2^{nd})$ grid points. Let $h_n = 2^{-n}$, therefore the sparse grid employs $O(h_n^{-1} \cdot \log_2(h_n^{-1})^{d-1})$ grid points compared to $O(h_n^{-d})$ nodes in the full grid.

Bungartz and Griebel \cite{bungartzGriebel} show that the accuracy of the sparse grid using $O(h_n^{-1} \cdot \log_2(h_n^{-1})^{d-1})$ nodes is of order $O(h_n^2 \log_2(h_n^{-1})^{d-1}))$ in the case of piecewise linear finite elements discretization and under the smoothness condition that the mixed derivatives are bounded. Thus, the accuracy of the sparse grid is only slightly deteriorated from the accuracy $O(h_n^{2})$ of conventional full grid methods which need $O(h_n^{-d})$ grid points. Therefore, sparse grids need much less points than regular full grids to achieve a similar approximation quality.

However, the structure of a sparse grid is more complicated than the one of a full grid. Common PDE solvers usually manage only full grid solutions. Existing sparse grid methods working directly in the hierarchical basis involve a challenging implementation \cite{sAchatz,aZeiser}. This handicap can be circumvented with the help of the sparse grid combination technique which not only exploits the economical structure of the sparse grids but also allows for the use of traditional full grid PDE solvers.

Concerning finite differences, in \cite{bgrz94, bgrz96} the authors obtain error bounds in terms of the Fourier transform coefficients when using a combination technique with a central difference scheme for the Laplace equation. The adaptive case with finite differences has been addressed in \cite{Griebel1998} for elliptic operators and the errors are obtained in $L^2$ and $L^{\infty}$-norms under the assumption that the $2d$-th mixed derivatives are bounded. In \cite{Koster}, the smoothness condition is written in terms of Holder spaces, the continuity of the mixed derivatives and their associated semi-norm to be finite, also including more general finite differences schemes. More recently, in \cite{reisinger} a methodology to obtain error bounds for general finite differences schemes in any dimension is proposed. It is mainly based on an error correction scheme leading to an appropriate error expansion. The results are again based on the existence of bounded mixed derivatives of the solution in the $L^{\infty}$ norm, which vanish at the boundary to avoid regularity problems. For this kind of functions, results about the expansion and interpolation are previously stated.

Concerning the smoothness of solution of the here treated PDEs, we note that the parabolic operator involves a smoothness effect (it is a kind of Black-Scholes operator in high spatial dimensions), so that the solution becomes $C^{\infty}$ for $t <T$, and the usual payoff (final condition at $t=T$) is continuous. Note that sparse grids have been analyzed when applied to basket options in \cite{resingerWittum}, for example.

\begin{figure}[!htb]
\centering
  \subfigure[$\Omega_s^{5}$, 177 grid points.] {\includegraphics[height=4.95cm]{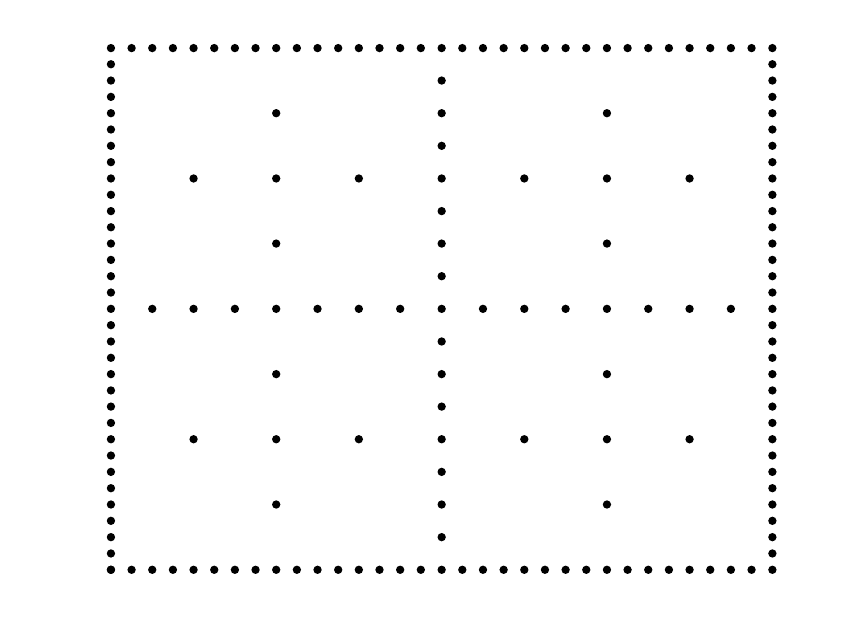}}
  \subfigure[$\Omega_s^{6}$, 385 grid points.] {\includegraphics[height=4.95cm]{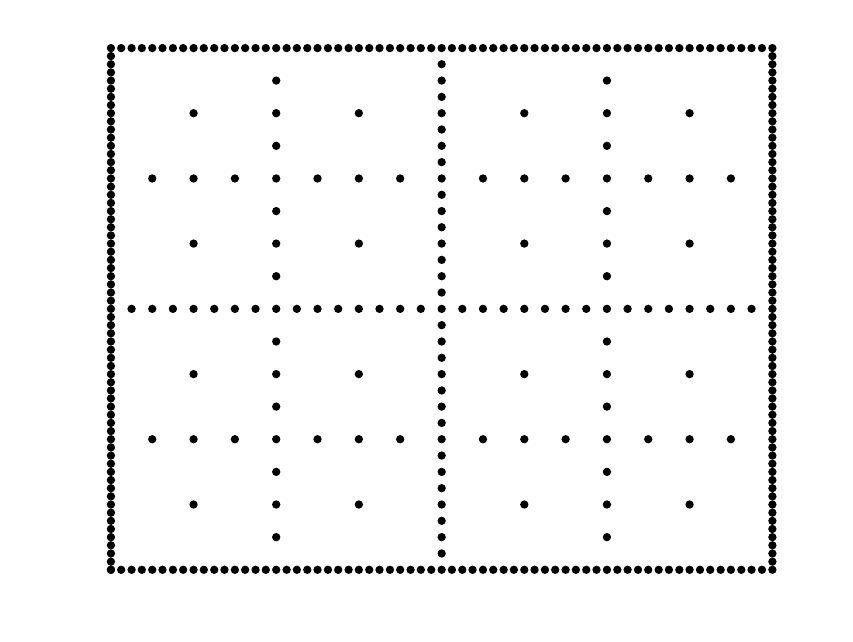}}
  \subfigure[$\Omega_s^{7}$, 833 grid points.] {\includegraphics[height=4.95cm]{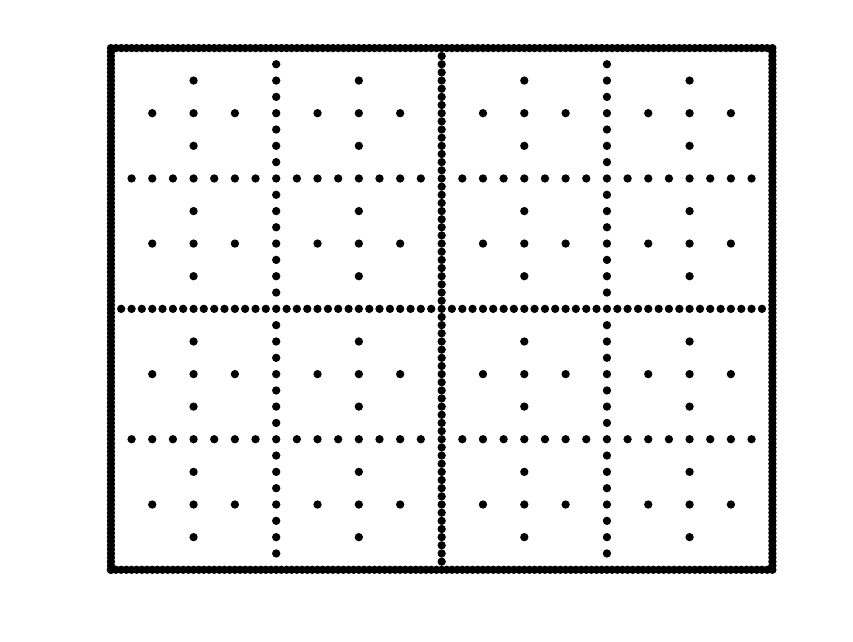}}
  \subfigure[$\Omega_s^{8}$, 1793 grid points.] {\includegraphics[height=4.95cm]{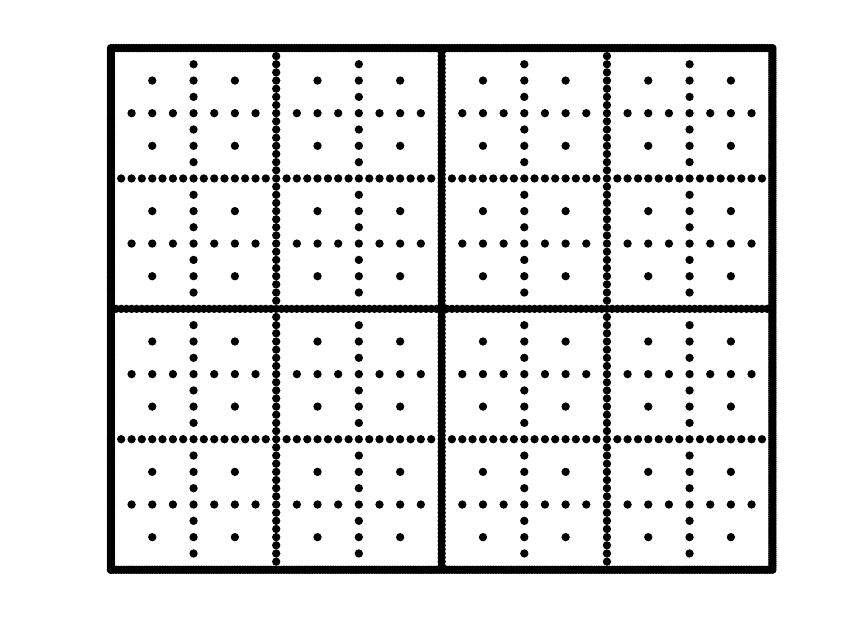}}
  \subfigure[$\Omega_s^{9}$, 3841 grid points.] {\includegraphics[height=4.95cm]{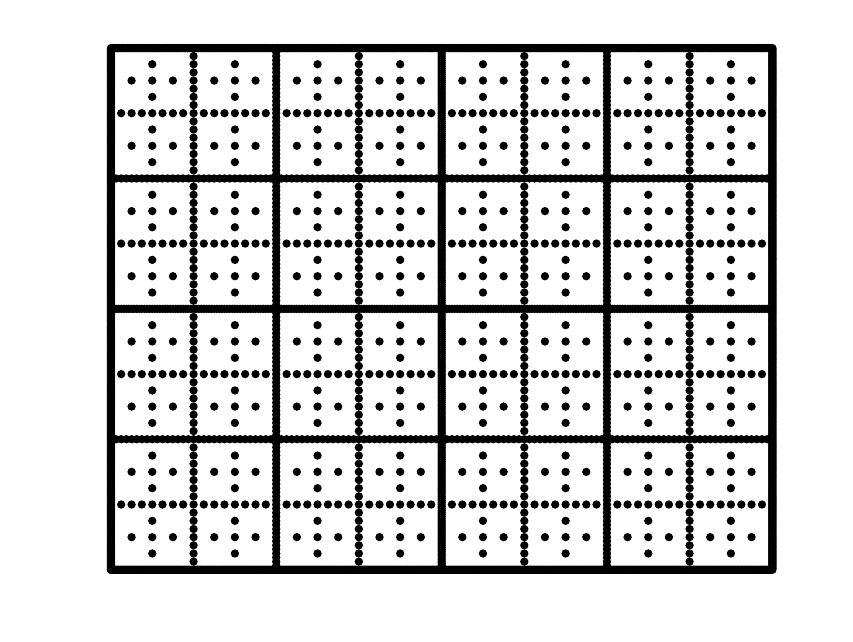}}
  \subfigure[$\Omega_s^{10}$, 8193 grid points.] {\includegraphics[height=4.95cm]{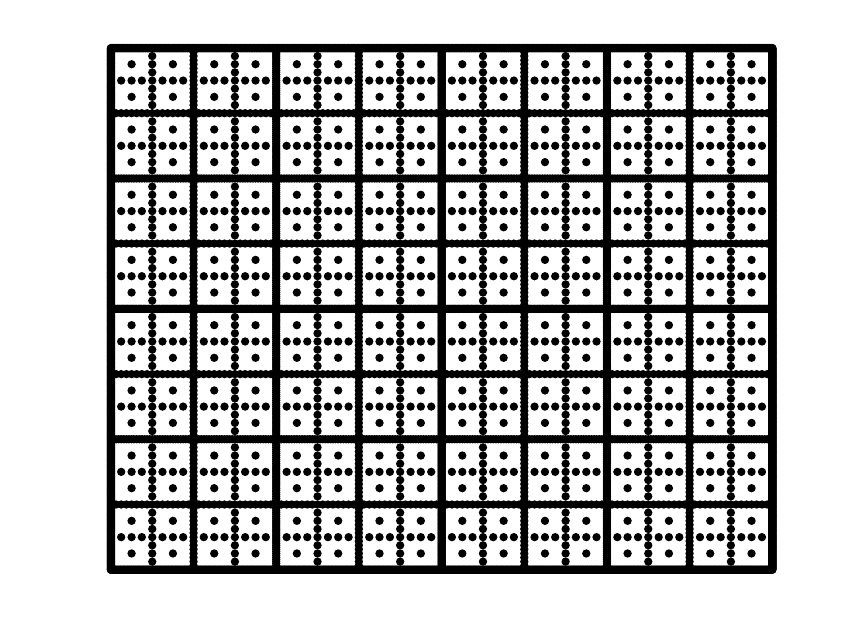}}
  \caption{Two dimensional sparse grids for levels $n=5,\ldots,10$.}
  \label{fig:merMorsparseGridD2}
\end{figure}

\begin{figure}[!htb]
\centering
  \subfigure[$\Omega_s^{5}$, 705 grid points.] {\includegraphics[height=5.10cm]{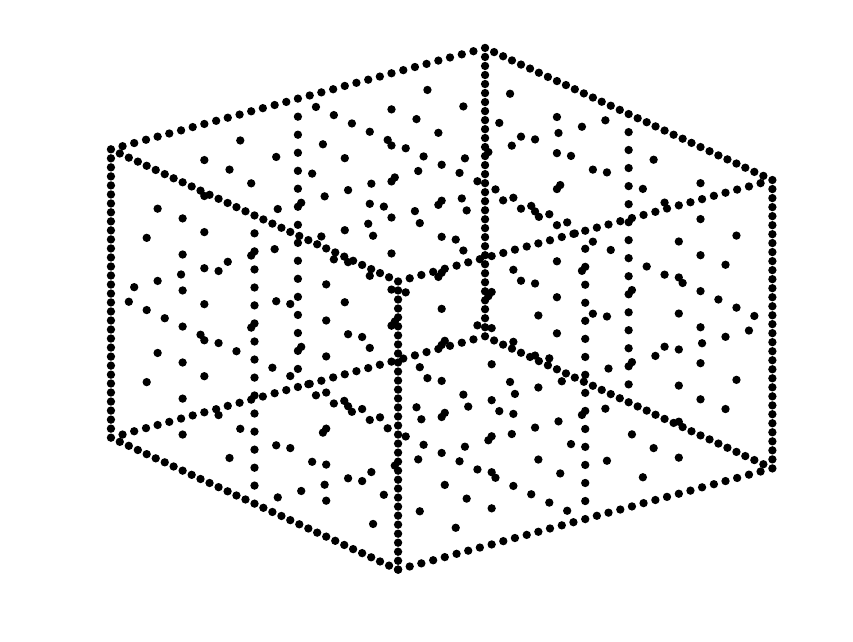}}
  \subfigure[$\Omega_s^{6}$, 1649 grid points.] {\includegraphics[height=5.10cm]{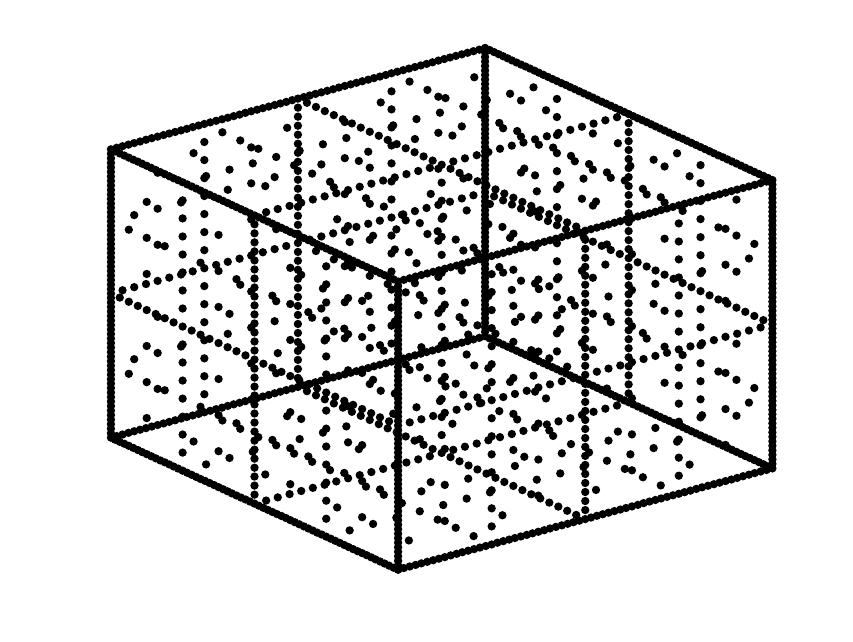}}
  \subfigure[$\Omega_s^{7}$, 3809 grid points.] {\includegraphics[height=5.10cm]{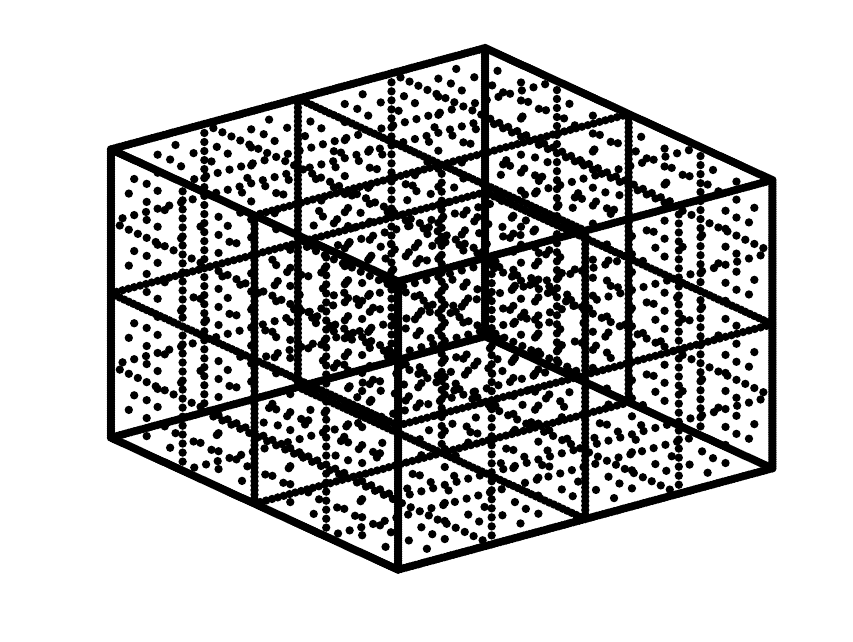}}
  \subfigure[$\Omega_s^{8}$, 8705 grid points.] {\includegraphics[height=5.10cm]{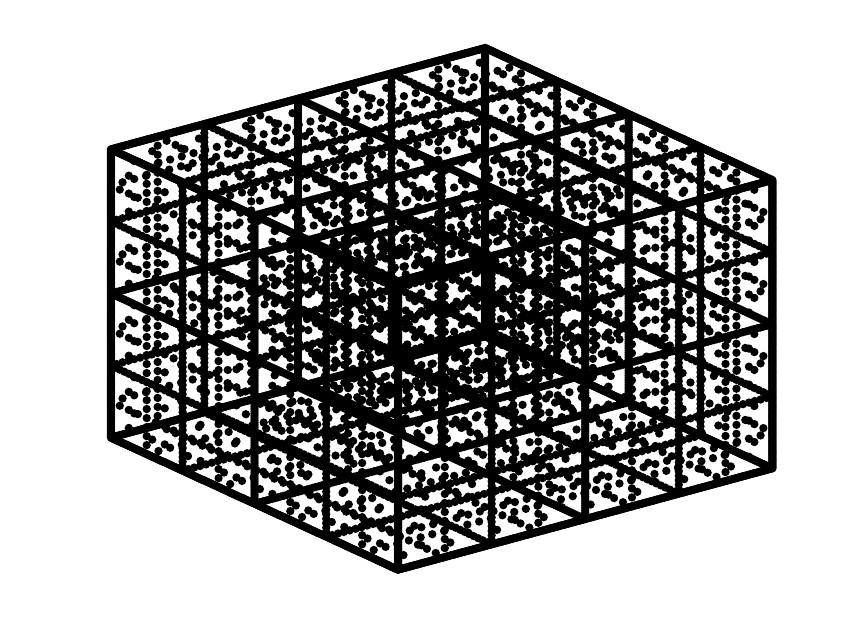}}
  \caption{Three dimensional sparse grids for levels $n=5, 6, 7$ and $8$.}
  \label{fig:merMorsparseGridD3}
\end{figure}

Finally, two and three dimensional sparse grids for several resolution levels $n$ are shown in Figures \ref{fig:merMorsparseGridD2} and \ref{fig:merMorsparseGridD3}, respectively. Additionally, the growth of the grid points when increasing $n$ can be observed.

\begin{figure}[!htb]
\begin{center}
\includegraphics[height=18cm]{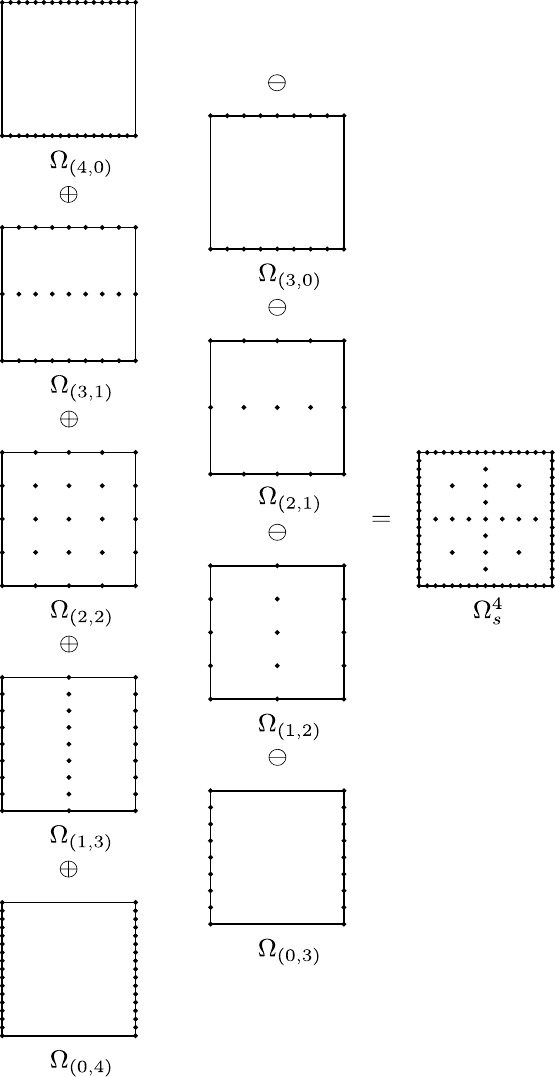}
\caption{Combination technique with level $n=4$ in two dimensions.}
\label{fig:merMorCombTech}
\end{center}
\end{figure}

\clearpage

\subsection{Combination technique}

Similar to the Richardson extrapolation \cite{richardson}, the so-called combination technique linearly combines the numerical solution on the sequence of anisotropic grids $\Omega_\mathbf{l}$ where
$$|\mathbf{l}|_1 = n-q, \quad q=0,\ldots,d-1.$$

The combination technique reads
\begin{equation}
 U_s^n = \displaystyle\sum_{q=0}^{d-1} (-1)^q \cdot {d-1 \choose q} \cdot \displaystyle \sum_{|\mathbf{l}|_1 = n-q} U_\mathbf{l}, \quad l_k\geq 0, \quad \forall k=1,\ldots,d,
\end{equation}
where $U_\mathbf{l}$ denotes the numerical solution on the grid $\Omega_\mathbf{l}$ and  $U_s^n$ the combined solution on the sparse grid $\Omega_s^n$.

The grids employed by the combination technique of level $n=4$ in two dimensions are shown in Figure \ref{fig:merMorCombTech}.

The idea of this technique is that the leading order errors from the dicretization on each grid cancel each other out in the combination solution.

The number of grid points involved in the approximation of $U_s^n$ grows according to $O(n^{d-1}\cdot 2^n)$. In fact, from the formula \eqref{eq:merMorNumOfSGP} we have to solve $n+d-1\choose d-1$ problems with $O(2^n)$ unknowns, $n+d-2 \choose d-1$ problems with $O(2^{n-1})$ unknowns, ... and $n \choose d-1$ problems with $O(2^{n-(d-1)})$ unknowns. This results in a total number of $O(n^{d-1}\cdot 2^n)$ grid points which is much less than the $O(2^{n\cdot d})$ grid nodes used by traditional full grid methods. Thus, the efficient use of sparse grids greatly reduces the computing time and the storage requirements which allows for the treatment of problems with ten variables and even more \cite{bungartzGriebel}.

We have seen that the combination technique linearly combines the numerical solution on several traditional full grids. The solution can be calculated on each of these grids by using any existing PDE numerical method like finite differences, finite volume or finite elements. In addition, since all these sub-problems are independent the combination technique can be parallelized \cite{griebel,heene}.

The combination technique approach presumes the existence of a so-called error splitting. It requires for an associated numerical approximation method on the full grid $\Omega_\mathbf{l}$ an error splitting of the form
\begin{equation}
 u(\mathbf{x}) - U_\mathbf{l}(\mathbf{x}) = \displaystyle\sum_{k=1}^d \sum_{\substack{\{j_1,\ldots,j_k\} \\ \subseteq \{1,\ldots,d\}}} C_{j_1,\ldots,j_k}(\mathbf{x},h_{{j_1}},\ldots,h_{{j_k}}) \cdot  h^p_{{j_1}} \cdot \ldots \cdot h^p_{{j_k}}, \label{errorSplitting}
\end{equation}
at each grid point $\mathbf{x} \in \Omega_\mathbf{l}$. Here $u$ denotes the exact solution of the partial differential equation under consideration, $U_\mathbf{l}$ the numerical solution on the grid $\Omega_\mathbf{l}$, $p>0$ is the order of accuracy of the numerical approximation method with respect to each coordinate direction and the coefficient functions $C_{j_1,\ldots,j_k}$ of $\mathbf{x}$ and the mesh sizes $h_{j_k}$, $k=1,\ldots,d$ are required to be bounded by a positive constant $K$ such that $$|C_{j_1,\ldots,j_k}(\mathbf{x},h_{{j_1}},\ldots,h_{{j_k}})| \leq K, \quad \forall k, 1\leq k \leq d, \quad\forall \{j_1,\ldots,j_m\} \subseteq \{1,\ldots,d\}. $$

The existence of the expansion (\ref{errorSplitting}) is a crucial point to obtain the error bounds of the sparse grid recombination technique and usually requires the assumption of bounded mixed derivatives.

In \cite{griebelThurner} Griebel and Thurner showed that if the solution of the PDE is sufficiently smooth, the pointwise accuracy of the sparse grid combination technique is $O(n^{d-1} \cdot 2^{-n\cdot p}) = O([\log_2h_n^{-1}]^{d-1} h_n^p)$, which is only slightly worse than $O(2^{-n\cdot p}) = O(h_n^p)$ obtained by the full grid solution.

The solution at points which do not belong to the sparse grid can be computed through interpolation. The applied interpolation method should provide at least the same order of accuracy of the numerical discretization scheme used to solve the PDE. Otherwise, the accuracy of the numerical scheme will be deteriorated.

Up to now we have assumed the existence of an error splitting of type \eqref{errorSplitting}. However, such an error splitting has to be proven for each problem. Nevertheless, proving the existence of this error splitting is usually very complex. Bungartz et al. \cite{bgrz94, bgrz96} showed the existence of such an error splitting for the finite differences discretization of the 2-d Laplace equation. Arciniega and Allen \cite{aa04} proved the existence of this error splitting for the fully implicit as well as the Crank-Nicolson discretization scheme of the European call option. More recently, Reisinger \cite{reisinger} showed that such a splitting also holds for a wider class of linear PDEs, for example convection-diffusion equations. The author gives general conditions which need to be fulfilled to ensure the existence of the desired splitting structure: sufficiently smooth initial data and compatible boundary data, the existence of bounded mixed derivatives, a consistent numerical scheme which provides a truncation error of the desired splitting structure and stability of the discretization scheme. As a summary, we can say that the deduction of the error splitting formula is very complex and was until now only performed for some reference problems. However, we will see in the following Section \ref{numericalResults} that the numerical results for the sparse grid combination technique are promising, even for more complex financial products.

\subsection{Numerical results} \label{numericalResults}

Taking advantage of the previously described sparse grid combination technique, in this section we are pricing the same interest rate derivatives that have been valued in the former Section \ref{numericaResutsDF} where traditional full grid finite differences methods were considered. In addition to those products, we are going to price interest rate derivatives with up to seven underlying LIBOR rates and their stochastic volatility, showing that the sparse grid combination technique is able to cope with the curse of dimensionality up to a certain extent. As in the previous Section \ref{numericaResutsDF}, we will use Crank-Nicolson scheme, we will consider the Gauss-Seidel iterative solver and the same boundary conditions as in Section \ref{boundaryConditions}. In the present case, we are interested in the evaluation of the solution at a single point which corresponds with the value of the forward rates at time zero (see Table \ref{table:marketData}) and $V(0)=1$. The numerical solution on each grid handled by the combination technique is interpolated at this point using multilinear interpolation and then added up with the appropriate weights.

The sparse grid combination technique has been implemented to run on multicore CPUs. The program was optimized and parallelized using OpenMP \cite{ref:openmp}. CPU times, measured in seconds, correspond to executions using $24$ threads, so as to take advantage of Intel Hyperthreading. The speedups of the parallel version with respect to the pure sequential code are around $16$. To the best of our knowledge, graphic processor units (GPUs) are not well-suited to parallelize the combination technique, due to the fact that the different grids employed by the combination technique involve memory accesses patterns totally different, therefore, it is not possible to access the device memory in a coalesced way \cite{NVIDIA-2011}, thus GPU global memory can not serve threads in parallel. In this scenario, the GPU code will be ill performing. In the work \cite{gaikwadToke} the authors take advantage of GPUs to parallelize the solver of each full grid considered by the combination technique. However, they do not parallelize the combination technique itself.

In Table \ref{table:d2_0} a $1\times 1$ European swaption is priced. The exact price of this derivative is $0.659096$ basis points, as discussed in Section \ref{numericaResutsDF}. These results are to be compared with those of Table \ref{table:df2_0}, where it can be seen how the computational times and the grid points employed by the sparse grid combination technique have been substantially reduced.

\begin{table}[!htb]
\begin{center}
{\footnotesize
\begin{tabular}{|r|r|r|r|r|}
\hline
Level & Solution & Error & Time & Grid points \\
\hline
\hline
 $3$ & $6.864576$ & $6.205480$ & $0.04$ & $37$\\ %
\hline
 $4$ & $2.207696$ & $1.548600$ & $0.04$ & $81$\\ %
\hline
 $5$ & $1.107670$ & $0.448573$ & $0.05$ & $177$\\ %
\hline
 $6$ & $0.788659$ & $0.129562$ & $0.05$ & $385$\\ %
\hline
 $7$ & $0.668489$ & $0.009393$ & $0.06$ & $833$\\ %
\hline
 $8$ & $0.662096$ & $0.002999$ & $0.12$ & $1793$\\ %
\hline
 $9$ & $0.659715$ & $0.000618$ & $0.54$ & $3841$\\ %
\hline
 $10$ & $0.659287$ & $0.000191$ & $2.68$ & $8193$\\ %
\hline
 $11$ & $0.659127$ & $0.000031$ & $16.79$ & $9217$\\ %
\hline
\end{tabular}
}
\end{center}
\caption{Convergence of the sparse grid finite differences solution in basis points for $1$ LIBOR and stochastic volatility, $\sigma=0$, $V(0)=1$, $\beta=1$, $256$ time steps. Exact solution, $0.659096$ basis points.}
\label{table:d2_0}
\end{table}

Next, in Table \ref{table:d2_1} a $1\times 1$ European swaption is priced considering stochastic volatility. These results are to be compared with those of Table \ref{table:df2_1}.

\begin{table}[!htb]
\begin{center}
{\footnotesize
\begin{tabular}{|r|r|r|}
\hline
Level & Solution & Time \\
\hline
\hline
 $3$ & $6.243$ & $0.03$  \\ %
\hline
 $4$ & $3.653$ & $0.04$  \\ %
\hline
 $5$ & $2.199$ & $0.04$  \\ %
\hline
 $6$ & $2.069$ & $0.08$  \\ %
\hline
 $7$ & $1.779$ & $0.08$  \\ %
\hline
 $8$ & $1.720$ & $0.39$  \\ %
\hline
 $9$ & $1.681$ & $1.13$  \\ %
\hline
 $10$ & $1.668$ & $6.91$  \\ %
\hline
 $11$ & $1.662$ & $43.03$  \\ %
\hline
\end{tabular}
}
\end{center}
\caption{Convergence of the sparse grid finite differences solution in basis points for $1$ LIBOR and stochastic volatility, $\sigma=0.3$, $\phi_i = 0.4$, $V(0)=1$, $\beta=1$, $256$ time steps. $95$\% confidence interval with Monte Carlo simulation using $10^7$ paths, $[1.652,1.672]$ in basis points.}

\label{table:d2_1}
\end{table}

In the following Tables \ref{table:d3} and \ref{table:d4}, the pricing of $1\times 2$ and $1\times 3$ European swaptions taking into account stochastic volatility is shown, as in Tables \ref{table:df3} and \ref{table:df4}, respectively. For the higher resolution levels, the full grid method became very slow, while the sparse grid combination technique results much faster. Note that the combination technique is able to price successfully the $1\times 3$ European swaption, this was not attainable in Table \ref{table:df3}.

\begin{table}[!htb]
\begin{center}
{\footnotesize
\begin{tabular}{|r|r|r|}
\hline
Level & Solution & Time \\
\hline
\hline
 $5$ & $9.560$ & $0.12$  \\ %
\hline
 $6$ & $7.416$ & $0.12$  \\ %
\hline
 $7$ & $5.896$ & $0.17$  \\ %
\hline
 $8$ & $5.318$ & $0.46$  \\ %
\hline
 $9$ & $5.007$ & $1.05$  \\ %
\hline
 $10$ & $4.826$ & $5.34$  \\ %
\hline
 $11$ & $4.833$ & $31.42$  \\ %
\hline
 $12$ & $4.820$ & $197.90$  \\ %
\hline
\end{tabular}
}
\end{center}
\caption{Convergence of the sparse grid finite differences solution in basis points for $2$ LIBORs and stochastic volatility, $\sigma=0.3$, $\phi_i = 0.4$, $V(0)=1$, $\beta=1$, $256$ time steps. $95$\% confidence interval with Monte Carlo simulation using $10^7$ paths, $[4.800,4.844]$ in basis points.}
\label{table:d3}
\end{table}

\begin{table}[!htb]
\begin{center}
{\footnotesize
\begin{tabular}{|r|r|r|}
\hline
Level & Solution & Time \\
\hline
\hline
 $7$ & $5.872$ & $0.34$  \\ %
\hline
 $8$ & $13.279$ & $0.88$  \\ %
\hline
 $9$ & $7.466$ & $2.29$  \\ %
\hline
 $10$ & $8.642$ & $8.37$  \\ %
\hline
 $11$ & $9.809$ & $29.74$  \\ %
\hline
 $12$ & $9.686$ & $156.75$  \\ %
\hline
 $13$ & $8.694$ & $895.75$  \\ %
\hline
 $14$ & $8.671$ & $5725.09$  \\ %
\hline
\end{tabular}
}
\end{center}
\vspace{-0.5cm}
\caption{Convergence of the sparse grid finite differences solution in basis points for $3$ LIBORs and stochastic volatility, $\sigma=0.3$, $\phi_i = 0.4$, $V(0)=1$, $\beta=1$, $256$ time steps. $95$\% confidence interval with Monte Carlo simulation using $10^7$ paths, $[8.635,8.700]$ in basis points.}
\label{table:d4}
\end{table}

\begin{table}[!htb]
\begin{center}
{\footnotesize
\begin{tabular}{|r|r|r|}
\hline
Level & Solution & Time \\
\hline
\hline
 $9$ & $12.51$ & $8.29$  \\ %
\hline
 $10$ & $10.35$ & $21.60$  \\ %
\hline
 $11$ & $15.75$ & $64.52$  \\ %
\hline
 $12$ & $16.30$ & $223.30$  \\ %
\hline
 $13$ & $9.95$ & $921.53$  \\ %
\hline
 $14$ & $13.03$ & $4504.83$  \\ %
\hline
$15$ & $13.24$ & $25980.70$  \\ %
\hline
\end{tabular}
}
\end{center}
\vspace{-0.5cm}
\caption{Convergence of the sparse grid finite differences solution in basis points for $4$ LIBORs and stochastic volatility, $\sigma=0.3$, $\phi_i = 0.4$, $V(0)=1$, $\beta=1$, $256$ time steps. $95$\% confidence interval with Monte Carlo simulation using $10^7$ paths, $[13.20,13.29]$ in basis points.}
\label{table:d5}
\end{table}

\begin{table}[!htb]
\begin{center}
{\footnotesize
\begin{tabular}{|r|r|r|}
\hline
Level & Solution & Time \\
\hline
\hline
 $11$ & $23.71$ & $305.46$  \\ %
\hline
 $12$ & $9.45$ & $855.16$  \\ %
\hline
 $13$ & $20.38$ & $2394.99$ \\ %
\hline
 $14$ & $18.82$ & $7584.20$ \\ %
\hline
 $15$ & $18.60$ & $29529.46$ \\ %
\hline
 $16$ & $18.56$ & $160027.61$ \\ %
\hline
\end{tabular}
}
\end{center}
\caption{Convergence of the sparse grid finite differences solution in basis points for $5$ LIBORs and stochastic volatility, $\sigma=0.3$, $\phi_i = 0.4$, $V(0)=1$, $\beta=1$, $256$ time steps. $95$\% confidence interval with Monte Carlo simulation using $10^7$ paths, $[18.51,18.61]$ in basis points.}
\label{table:d6}
\end{table}

\begin{table}[!htb]
\begin{center}
{\footnotesize
\begin{tabular}{|r|r|r|}
\hline
Level & Solution & Time \\
\hline
\hline
 $13$ & $18.30$ & $627.97$  \\ %
\hline
 $14$ & $15.55$ & $2115.55$  \\ %
\hline
 $15$ & $25.60$ & $9892.70$  \\ %
\hline
 $16$ & $24.49$ & $50139.95$  \\ %
\hline
\end{tabular}
}
\end{center}
\caption{Convergence of the sparse grid finite diffferences solution in basis points for $6$ LIBORs and stochastic volatility, $\sigma=0.3$, $\phi_i = 0.4$, $V(0)=1$, $\beta=1$, $12$ time steps. $95$\% confidence interval with Monte Carlo simulation using $10^7$ paths, $[24.47,24.59]$ in basis points.}
\label{table:d7}
\end{table}

\begin{table}[!htb]
\begin{center}
{\footnotesize
\begin{tabular}{|r|r|r|}
\hline
Level & Solution & Time \\
\hline
\hline
 $15$ & $19.56$ & $24612.73$  \\ %
\hline
 $16$ & $26.80$ & $79463.11$  \\ %
\hline
 $17$ & $30.93$ & $324996.63$  \\ %
\hline
\end{tabular}
}
\end{center}
\caption{Convergence of the sparse grid finite differences solution in basis points for $7$ LIBORs and stochastic volatility, $\sigma=0.3$, $\phi_i = 0.4$, $V(0)=1$, $\beta=1$, $2$ time steps. $95$\% confidence interval with Monte Carlo simulation using $10^7$ paths, $[30.85,30.98]$ in basis points.}
\label{table:d8}
\end{table}


Finally, in Tables from \ref{table:d5} to \ref{table:d8}, $1 \times 4$, ..., $1 \times 7$ European swaptions are priced considering stochastic volatility. The pricing of these interest rate derivatives was not viable with the full grid approach of Section \ref{finiteDifferenceMethod}. In order to be able to price derivatives with more than $8$ underlyings, the combination technique method should be parallelized to run on a cluster of processors. In the Chapter 13 of the book \cite{gerstner}  Philipp Schröder et al. discuss the parallelization of the combination technique using MPI (\textit{Message Passing Interface}) API. In \cite{larson} the authors parallelize the sparse grid combination technique taking advantage of a MapReduce framework, algorithms that are inherently fault tolerant.

\section{Conclusion}
In this work we have posed for first time in the literature the PDEs associated to the SABR/LIBOR market models proposed by Mercurio \& Morini \cite{mercurioMorini} and Hagan \cite{haganSABRLIBOR}. In order to price interest rate derivatives we have developed a traditional full grid finite difference method. This approach is able to successfully price derivatives up to two or three underlying forward rates in reasonable computational times. However, when the number of underlyings increases this scheme suffers from the well-known curse of dimensionality. In order to price derivatives over a moderately large number of forward rates we have proposed to use the sparse grid combination technique. Taking into account that this technique is embarrassingly parallel  we have parallelized it so as to drastically reduce computational times. Finally, we have tested the proposed method in two different ways. On one hand, using the analytical solution when it exists. On the other hand, when the exact solution is not known, we have used as reference solution the one computed with the Monte Carlo method, thus ensuring the correctness of the developed scheme.




 \clearpage
\bibliographystyle{plain}
\bibliography{mybibfile}

\end{document}